\input ./style/arxiv-general.cfg
\documentclass[MSNbibl,citesort,number,dvips]{arxbj}
\makeatletter
   \@ifpackageloaded{graphicx}{}{\usepackage{graphicx}}
\makeatother
\usepackage{mathrsfs,upgreek}

\aid{0}
\volume{22}
\issue{1}
\pubyear{2016}
\firstpage{1}
\lastpage{37}
\doi{10.3150/14-BEJ609} 

\makeatletter
\mathchardef\varSigma="0106
\newcommand{\rrVert}{\Vert}
\newcommand{\llVert}{\Vert}

\newcommand{\rright}{\right}
\newcommand{\lleft}{\left}
\newcommand{\rrvert}{\vert}
\newcommand{\llvert}{\vert}
\newtheorem{cor}{Corollary}
\newtheorem{lemma}{Lemma}[section]
\newproclaim{defn}{Definition}
\newtheorem{prop}{Proposition}

\def\iidsim{\stackrel{\mathrm{i.i.d.}}{\sim}}
\def\Cov{\operatorname{Cov}}
\def\iidsim{\stackrel{\mathrm{i.i.d.}}{\sim}}
\def\O{\mathcal{O}}
\def\E{\mathscr{E}}
\def\z{\mathbf{z}}
\def\dd{\boldsymbol{\delta}}
\def\bth{\boldsymbol{\theta}}
\def\th{\theta}
\def\vv{\mathbf{v}}
\def\bxi{\boldsymbol{\xi}}
\def\tr{\operatorname{tr}}
\def\1{\mathbf{1}}
\def\ee{\boldsymbol{\epsilon}}
\def\bb{\boldsymbol{\beta}}
\def\d{\delta}
\def\Var{\operatorname{Var}}
\def\g{\gamma}
\def\x{\mathbf{x}}
\def\w{\mathbf{w}}
\def\F{\mathbb{F}}
\def\t{\mathbf{t}}
\def\u{\mathbf{u}}
\def\y{\mathbf{y}}
\def\t{\mathbf{t}}
\def\R{\mathbb{R}}
\def\S{\varSigma}

\renewcommand{\epsilon}{\varepsilon}
\def\sfrac#1#2{#1/#2}

\def\afrac#1#2{#1/(#2)}

\def\sklfrac#1#2{(#1/#2)}

\makeatother

\begin{document}
\begin{frontmatter}

\title{Ridge regression and asymptotic minimax~estimation over spheres of growing~dimension}
\runtitle{Ridge regression and asymptotic minimax estimation}

\begin{aug}
\author[A]{\inits{L.H.}\fnms{Lee H.} \snm{Dicker}\corref{}\ead[label=e1]{ldicker@stat.rutgers.edu}}
\address[A]{Department of Statistics and Biostatistics, Rutgers
University, 501 Hill Center,
110 Frelinghuysen Road, Piscataway, NJ 08854, USA.
\printead{e1}}
\end{aug}

\received{\smonth{1} \syear{2013}}
\revised{\smonth{9} \syear{2013}}

%
\begin{abstract}
We study asymptotic minimax problems for estimating a $d$-dimensional
regression parameter over spheres of growing dimension ($d \to
\infty$). Assuming that the data follows a linear model with Gaussian
predictors and errors, we show that ridge regression is asymptotically
minimax and
derive new closed form expressions for its
asymptotic risk under squared-error loss. The asymptotic risk of
ridge regression is closely related to the Stieltjes transform of the
Mar\v{c}enko--Pastur distribution and the spectral distribution of the
predictors from the linear model. Adaptive ridge
estimators are also proposed (which adapt to the unknown radius of the
sphere) and connections with equivariant
estimation are highlighted. Our results are
mostly relevant for asymptotic settings where the number of
observations, $n$, is proportional to the number of predictors,
that is, $d/n \to\rho\in(0,\infty)$.
\end{abstract}

%

\end{frontmatter}

\section{Introduction}\label{sec1}

Consider a linear model where the observed data consists of outcomes
$y_1,\ldots,y_n\in\R$ and
\mbox{$d$-dimensional} predictors $\x_1,\ldots,\x_n \in\R^d$ that are related
via the equation
%
%
\begin{equation}
\label{lm} y_i = \x_i^T\bb+ \epsilon
_i,\qquad  i = 1,\ldots,n;
\end{equation}
the $d$-dimensional vector $\bb= (\beta_1,\ldots,\beta_d)^T \in\R^d$
is an unknown parameter
and $\epsilon_1,\ldots,\epsilon_n\in\R$ are
unobserved errors. To simplify notation, let $\y= (y_1,\ldots,y_n)^T
\in\R^n$, $X = (\x_1,\ldots,\x_n)^T$, and $\ee=
(\epsilon_1,\ldots,\epsilon_n)^T\in\R^n$. Then (\ref{lm}) may be
rewritten as $\y= X\bb+ \ee$.

In this paper, we study asymptotic minimax estimation of $\bb$ over
spheres of growing dimension (i.e., $d \to\infty$), under the
assumption that the data $(\y,X)$ are jointly Gaussian. This is a
variant of a problem considered
by Goldenshluger and Tsybakov \cite
{goldenshluger2001adaptive,goldenshluger2003optimal}; it is
closely related to the fundamental work of
Pinsker \cite{pinsker1980optimal} and others,
for example, Belitser and Levit \cite{belitser1995minimax}, Beran \cite{beran1996stein},
Golubev \cite{golubev1987adaptive},
on sharp
asymptotic minimax estimation in the Gaussian
sequence model. Taken together, the results in this paper provide a new example where sharp
asymptotic minimax estimation is possible; an example that illustrates
connections between linear models with many predictors and now
classical results on the spectral distribution
of large random matrices.

\subsection{Statement of problem}\label{sec1.1}

Let $I_k$ denote the $k \times
k$ identity matrix. We assume throughout that
%
%
\begin{equation}
\label{norm} \x_1,\ldots,\x_n \iidsim
N(0,I_d)\quad  \mbox{and}\quad  \ee\sim N(0,I_n)
\end{equation}
are independent. More general models, where one might allow for positive
definite $\Cov(\x_i) = \S$ and arbitrary $\Var(\epsilon_i) =
\sigma^2 > 0$, are
discussed in Section~\ref{sec1.4}.

Given an estimator
$\hat{\bb} = \hat{\bb}(\y,X)$ of $\bb$, define the
risk under squared-error loss
%
%
\begin{equation}
\label{risk} R(\hat{\bb},\bb) = E_{\bb} \bigl(\llVert \hat{\bb} - \bb
\rrVert ^2 \bigr),
\end{equation}
where $\llVert  \cdot\rrVert  $ denotes the $\ell^2$-norm. The expectation in
(\ref{risk}) is taken with respect to the joint distribution of
$(X,\ee)$
and the subscript $\bb$
in $E_{\bb}$ indicates that $\y= X\bb+ \ee$ (for expectations that
do not involve $\y$, we will often omit the subscript). We
emphasize that the expectation in (\ref{risk}) is taken over the
predictors $X$ as well as the errors $\ee$; in other words, rather
than conditioning on $X$,
(\ref{risk})~is the \emph{unconditional} risk under squared-error
loss.

Let $S^{d-1}(\tau) = \{\bb\in\R^d; \llVert  \bb\rrVert   = \tau\}$ be the
sphere of radius $\tau\geq0$ in $\R^d$ centered at the origin. The
minimax risk for estimating $\bb$
over $S^{d-1}(\tau)$ is
given by
%
%
\begin{equation}
\label{mr} r(\tau) = r_{d,n}(\tau) = \inf_{\hat{\bb}}
\sup_{\bb\in
S^{d-1}(\tau)} R(\hat{\bb},\bb),
\end{equation}
where the infimum
in (\ref{mr}) is taken over all measurable estimators $\hat{\bb} =
\hat{\bb}(\y,X)$.

The minimax problem
determined by (\ref{mr}) is the main focus of this paper. Our analysis
entails (i)
identifying and analyzing specific
estimators $\hat{\bb}$ such that $\sup_{\bb\in S^{d-1}(\tau)}
R(\hat{\bb},\bb) \approx r(\tau)$, and (ii) obtaining accurate closed-form
approximations for $r(\tau)$, while focusing on settings where $d$ is large.

\subsection{Overview of results}\label{sec1.2}

To better orient the reader, we give a brief
section-by-section overview of the paper. We conclude this section with an
additional comment on the nature of the asymptotic results derived
herein.

\emph{Section~\ref{sec2}}: \emph{Ridge regression}. Ridge regression
(Hoerl and Kennard \cite{hoerl1970ridge} and Tihonov \cite{tikhonov1963solution}) is a widely studied
regularized estimation method whose use has been advocated in various settings
where $d$ is large or $X$ is ill-conditioned. Our analysis in Section~\ref{sec2}
yields a simple formula for the optimal ridge
regularization parameter and a new closed-form expression for the
associated ridge estimator's asymptotic risk. More specifically, we
show that if $d/n \to\rho\in(0,\infty)$, then the asymptotic risk
of the ridge estimator is closely related to the Stieltjes transform of the
Mar\v{c}enko--Pastur distribution (Mar{\v{c}}enko and Pastur \cite{marcenko1967distribution}),
which plays a prominent role in random matrix theory, for example,
Bai \textit{et al.} \cite{bai2003convergence},
El~Karoui \cite{elkaroui2008spectrum},
Silverstein \cite{silverstein1995strong}.
Settings where
$d/n \to0$ and $d/n \to\infty$ are also considered. Our results
for ridge regression immediately provide an upper bound on $r(\tau)$, in
the usual way: It is clear from (\ref{mr}) that $r(\tau)
\leq\sup_{\bb\in S^{d-1}(\tau)} R(\hat{\bb},\bb)$
for all estimators $\hat{\bb}$; taking $\hat{\bb}$ to be the specified
ridge estimator gives the desired upper bound.

\emph{Section~\ref{sec3}}: \textit{An equivalent Bayes problem}. An equivariance
argument implies that $r(\tau)$ is
equal to the Bayes risk for estimating $\bb$ under the prior distribution
$\bb\sim\pi_{S^{d-1}(\tau)}$, where $\pi_{S^{d-1}(\tau)}$
denotes the uniform distribution on $S^{d-1}(\tau)$ (this is an application
of well-known results on equivariance, e.g., Chapter~6 of
Berger \cite{berger1985statistical}, and
is essentially an
illustration of the
Hunt--Stein theorem (Bondar and Milnes \cite
{bondar1981amenability})). Additionally, we
argue that when $d$ is large, the Bayes risk for estimating
$\bb$ under the prior distribution $\bb\sim
\pi_{S^{d-1}(\tau)}$ is close to the Bayes risk for estimating
$\bb$ under a normal prior distribution, which coincides with the risk
of ridge regression. We conclude that the risk of ridge regression is
asymptotically equivalent to $r(\tau)$ and that ridge
regression is asymptotically optimal for estimation over $S^{d-1}(\tau)$.

\emph{Section~\ref{sec4}}: \textit{An adaptive ridge estimator}. The ridge regression
estimator $\hat{\bb}_r(\tau)$ that is asymptotically optimal over
$S^{d-1}(\tau)$ depends on the radius
$\tau= \llVert  \bb\rrVert  $, which is typically unknown. Replacing $\tau$ with
an estimate, we obtain
an adaptive ridge estimator that does not depend on $\tau$, but is
asymptotically equivalent to $\hat{\bb}_r(\tau)$. It follows that
the adaptive ridge
estimator is adaptive asymptotic minimax over spheres $S^{d-1}(\tau)$, provided
$\tau^2 \gg n^{-1/2}$. Additionally, we show that the
adaptive ridge estimator is asymptotically optimal among the class of
all estimators for $\bb$ that are equivariant with respect to
orthogonal transformations of the predictors, as $d \to\infty$.

Proofs may be found in the \hyperref[appA]{Appendices}.

\emph{Note on asymptotics}. Throughout the paper, our asymptotic
analysis is focused on settings where $d \to\infty$. We typically
assume that $n \to\infty$ along with $d$ and that $d/n \to\rho\in
[0,\infty]$. It will become apparent
below that most of the ``action'' occurs when $0 < \rho< \infty$.
Indeed, one of the implications of our results is that if $0 <
\rho< \infty$, then the minimax risk $r(\tau)$ is influenced by the
spectral distribution of the empirical
covariance matrix $n^{-1}X^TX$. On the other hand, if $\rho= 0$, then
the behavior of $r(\tau)$ is more
standard. If $\rho= \infty$, then we will show that it is impossible
to out-perform the
trivial estimator $\hat{\bb}_{\mathrm{null}} = 0$ for estimation over
$S^{d-1}(\tau)$; note the contrast with sparse estimation problems,
where $\bb$ is assumed to be sparse and it may be possible to
dramatically out-perform $\hat{\bb}_{\mathrm{null}}$ when $d/n \to\infty$,
for example, Bickel \textit{et al.} \cite{bickel2009simultaneous},
Bunea \textit{et al.} \cite{bunea2007sparsity}, Candes and Tao \cite{candes2007dantzig},
Raskutti \textit{et al.}
\cite{raskutti2011minimax}, Ye and Zhang \cite{ye2010rate},  Zhang \cite{zhang2010nearly}.

\subsection{Relationship to existing work}\label{sec1.3}

The minimax problem (\ref{mr}) is closely
related to problems considered by Goldenshluger and Tsybakov \cite{goldenshluger2001adaptive,goldenshluger2003optimal}, who
studied minimax prediction problems over $\ell^2$-ellipsoids
\[
\Biggl\{\bb\in\ell^2; \sum_{k = 1}^{\infty}
a_k^2\beta_k^2 \leq
L^2 \Biggr\};\qquad  L > 0, a = \{a_k\},
\]
in an infinite-dimensional linear model with independent (but not
necessarily Gaussian) predictors. Goldenshluger and Tsykbakov's results
apply to classes of
ellipsoids with various constraints on $a = \{a_k\}$ and $L$.
Taking $L =\tau$, $a_1 = \cdots= a_d = 1$, and $a_{d+1} = a_{d+2} =
\cdots=
\infty$ (and following the convention that $0 \times\infty= 0$), the
results in Goldenshluger and Tsybakov \cite
{goldenshluger2003optimal} may be
applied to obtain asymptotics for the
minimax risk over the $d$-dimensional ball $B_d(\tau) = \{\bb\in\R
^d; \llVert  \bb\rrVert   \leq
\tau\}$,
%
%
\begin{equation}
\label{mrb} \bar{r}(\tau) = \inf_{\hat{\bb}} \sup
_{\bb\in B_d(\tau)} R(\hat{\bb},\bb);
\end{equation}
results in Goldenshluger and Tsybakov \cite
{goldenshluger2001adaptive} yield adaptive
estimators that are asymptotically minimax over classes of balls
$B_d(\tau)$.
In Section~\ref{sec3.2}, we show that $r(\tau) \approx\bar{r}(\tau)$, when
$d$ is
large (see (\ref{mrbmr}) below). Thus, Goldenshluger and Tsybakov's
results are clearly related to the results presented here. However,
as applied to balls $B_d(\tau)$, their results typically require that
$d/n \to0$ (for
instance, Theorem~1 of Goldenshluger and Tsybakov \cite
{goldenshluger2003optimal}
requires that $d = \mathrm{o}\{\sqrt{n/\log(n)}\}$ and Assumption~3 of
Goldenshluger and Tsybakov \cite
{goldenshluger2001adaptive} requires $d = \mathrm{O}\{\sqrt{n}/\log(n)
\}$). By contrast, the
results in this paper apply in settings where $d/n \to\rho\in
[0,\infty]$, with the bulk of our work focusing on $0 < \rho<
\infty$.

The analysis in this paper focuses on estimation
over the sphere $S^{d-1}(\tau)$, rather than the ball $B_d(\tau)$;
that is, we
focus on the minimax problem (\ref{mr}), as opposed to (\ref{mrb}).
The ball
$B_d(\tau)$ and other star-shaped parameter spaces (e.g., ellipsoids or
$\ell^p$-balls) have
been more frequently studied in the literature on asymptotic minimax
problems over restricted parameter spaces
(Donoho and Johnstone \cite{donoho1994minimax},
Golubev \cite{golubev1990quasilinear},
Nussbaum \cite{nussbaum1999minimax}).
Evidently, the problems (\ref{mr})
and (\ref{mrb}) are closely related. However, analysis of (\ref{mr})
appears to be somewhat more complex; in particular, obtaining lower
bounds on $r(\tau)$ seems more challenging. To justify our emphasis on the
sphere $S^{d-1}(\tau)$, in Section~\ref{sec3.2}, we show that asymptotics for
$\bar{r}(\tau)$
follow easily from asymptotics for $r(\tau)$. Additionally, by
studying estimation over the sphere, we are able to draw deeper
connections with equivariance than seem to be
available if one focuses on the ball (e.g., Proposition~\ref{Oadapt}
below). A~similar approach has been considered by
Marchand \cite{marchand1993estimation} and
Beran \cite{beran1996stein} in their
analysis of the finite-dimensional Gaussian sequence mode. In fact,
one of the key technical results in this paper (Theorem~\ref{main}) is
essentially a multivariate extension of Theorem~3.1 in
Marchand \cite{marchand1993estimation}.
While we
believe that the additional insights provided by studying minimax
problems over the sphere justify the added complexity, we also note
that more standard approaches to obtaining lower bounds on the minimax
risk over
balls (see, e.g., Nussbaum \cite
{nussbaum1999minimax} or Chapter~3 of
Tsybakov \cite
{tsybakov2009introduction}) may be applied to obtain lower
bounds for
$\bar{r}(\tau)$ directly.

Finally in this section, we mention some of the existing work on
random matrix theory that is especially relevant for our analysis of
ridge regression in Section~\ref{sec2}. Theorem~\ref{mprisk} in Section~\ref{sec2.2}
relies heavily on now classical
results that describe the asymptotic
behavior of the empirical
distribution of the eigenvalues of $n^{-1}X^TX$ in high dimensions
(Bai \cite{bai1993convergence}, Bai \textit{et al.} \cite{bai2003convergence},
Mar{\v{c}}enko and Pastur \cite{marcenko1967distribution}).
Additionally, we point out that while other
authors have alluded to the relevance of random matrix
theory for ridge regression (El Karoui and K{\"o}sters \cite
{elkaroui2011geometric}), the
results presented here on ridge regression's asymptotic risk seem to
provide a greater level of detail than available elsewhere, in the
specified setting.

\subsection{Distributional assumptions}\label{sec1.4}

The linear model (\ref{lm}) with distributional assumptions (\ref
{norm}) is highly specialized. However, similar models have been
studied previously. Stein \cite
{stein1960multiple}, Baranchik \cite
{baranchik1973inadmissibility},
Breiman and Freedman \cite{breiman1983many}, Brown \cite{brown1990ancillarity} and Leeb \cite{leeb2009conditional}
studied estimation problems for linear models with jointly Gaussian
data, but, for the most part, these authors do not require $\Cov(\x_i)
= I_d$. Moreover, as discussed in Section~\ref{sec1.3}, the infinite-dimensional linear
model considered by
Goldenshluger and Tsybakov \cite
{goldenshluger2001adaptive,goldenshluger2003optimal} is
similar to the model studied in this paper. For our purposes, one of
the more significant consequences of the normality assumption
(\ref{norm}) is that the distributions of $X$ and $\ee$ are invariant
under orthogonal transformations. This leads to substantial
simplifications in many of the ensuing calculations. Results in El Karoui and
K{\"o}sters \cite{elkaroui2011geometric}
suggest that a general approach to
relaxing some of the distributional assumptions made in this paper may
be feasible, but this is not pursued further here.

We point out that the assumption $E(\x_i) = 0$, which is implicit in
(\ref{norm}), is not particularly limiting:
If $E(\x_i) \neq0$, then we can reduce to the mean 0 case by
centering and de-correlating the data. The normality assumption
(\ref{norm}) also requires $\Var(\epsilon_i) = 1$. If $\Var
(\epsilon_i) = \sigma^2 \neq1$
and $\sigma^2$ is known, then this can be reduced to the case where
$\Var(\epsilon_i) = 1$ by transforming the data $(\y,X) \mapsto(\y
/\sigma,X)$; the corresponding transformation for the parameters $\bb
$, $\sigma^2$ is
given by $(\bb,\sigma^2) \mapsto(\bb/\sigma,1)$ and the risk
function should be scaled by $\sigma^2$,
as well (ultimately in this scenario, most of the results in this
paper remain valid except that the signal-to-noise ratio
$\llVert  \bb\rrVert  ^2/\sigma^2$ replaces the signal strength $\llVert  \bb\rrVert  ^2$). If
$\sigma^2$ is unknown and $d/n \to\rho< 1$, then
$\sigma^2$ may be
effectively estimated by $\hat{\sigma}^2 = (n-d)^{-1}\llVert  y -
X\hat{\bb}_{\mathrm{ols}}\rrVert  ^2$, where $\hat{\bb}_{\mathrm{ols}} = (X^TX)^{-1}X^T\y$ is
the ordinary least squares (OLS) estimator; one can subsequently reduce
to the case where $\Var(\epsilon_i)
= 1$. (Throughout, if the square matrix $A$ is not invertible, then we
take $A^{-1}$ to be
its Moore--Penrose pseudoinverse; typically,
the matrices we seek to invert will be invertible with probability 1.)
Recent work suggests that $\sigma^2$ may also be effectively estimated
when $d > n$. Fan \textit{et al.} \cite{fan2012variance} and
Sun and Zhang \cite{sun2012scaled} propose
methods for estimating $\sigma^2$ when $d > n$ and $\bb$ is sparse (see
also related work by Belloni \textit{et al.} \cite
{belloni2011square} and
Dalalyan and Chen \cite{dalalyan2012fused} on estimating
$\bb$ in high dimensions when
$\sigma^2$ is unknown); Dicker \cite
{dicker2012residual}
considers estimating $\sigma^2$ when $d >n$ and $\bb$ is not sparse.

Under the Gaussian assumption (\ref{norm}), the predictors $\x_i$ are
uncorrelated at the population level, that is, $\Cov(\x_i) = I_d$. The
results in
this paper are easily adapted to settings
where $\Cov(\x_i) = \S$ is a known positive definite matrix by
transforming the data $(\y,X)
\mapsto(\y,X\S^{-1/2})$, and making corresponding transformations of
the parameters and risk function. If $\Cov(\x_i) = \S$ is unknown,
but $\hat{\S}$ is an operator norm
consistent estimator, then it is straigthforward to check that most of
our asymptotic results
remain valid, mutatis mutandis, for the transformed data
$(\y,X\hat{\S}^{-1/2})$. On the other hand, in high-dimensional
settings where $d/n \to
\rho> 0$, an operator norm consistent estimator for $\S$ may not exist.
In Dicker \cite{dicker2012optimal}, the
author considers a prediction
problem closely related to the estimation
problem considered in this paper, with unknown $\Cov(\x_i) = \S$;
the author identifies an asymptotically optimal
equivariant estimator and derives expressions for the estimator's
asymptotic risk (Theorems 2--3 and Corollary~1 of Dicker \cite{dicker2012optimal}). One interpretation of
the results in
Dicker \cite{dicker2012optimal} is that they
quantify the loss in efficiency of equivariant estimators when
$\Cov(\x_i) = \S$ is unknown, as compared to the results presented
here for the case where $\Cov(\x_i) = I_d$ is known.

\section{Ridge regression}\label{sec2}

Define the ridge regression estimator
\[
\hat{\bb}_r(t) = \bigl(X^TX + d/t^2
I_d \bigr)^{-1}X^T\y,\qquad  t \in[0,\infty].
\]
The parameter $t$ is referred to as the ``regularization'' or
``ridge'' parameter and is subject to further specification. By
convention, we take $\hat{\bb}_r(0) = 0$ and $\hat{\bb}_r(\infty)
= \hat{\bb}_{\mathrm{ols}} =
(X^TX)^{-1}X^T\y$ to be the OLS estimator.

\subsection{The oracle ridge estimator}\label{sec2.1}

Our first result identifies the optimal
ridge parameter $t$ and yields an oracle ridge
estimator with minimal risk. A simplified expression for the oracle
ridge estimator's risk is also provided.

%
\begin{prop}\label{rrisk}
Suppose that $\bb\in S^{d-1}(\tau)$. Then
%
%
\begin{equation}
\label{rrisk1} R \bigl\{\hat{\bb}_r(\tau),\bb \bigr\} = \inf
_{t \in[0,\infty]} R \bigl\{\hat{\bb}_r(t),\bb \bigr\} = E
\bigl[ \tr \bigl\{ \bigl(X^TX + d/\tau^2I_d
\bigr)^{-1} \bigr\} \bigr].
\end{equation}
\end{prop}

%
\begin{cor}\label{ubcor} Suppose that $\tau\geq0$. Then
\[
r(\tau) \leq\sup_{\bb\in S^{d-1}(\tau)} R \bigl\{\hat{\bb }_r(\tau
),\bb \bigr\} = E \bigl[\tr \bigl\{ \bigl(X^TX + d/
\tau^2I_d \bigr)^{-1} \bigr\} \bigr].
\]
\end{cor}

Proposition~\ref{rrisk} is proved in Appendix \hyperref[appA]{A} and it implies that
the optimal ridge parameter is given by the signal strength
$\tau= \llVert  \bb\rrVert  $. Notice that the risk of
$\hat{\bb}_r(\tau)$ is
constant over the sphere $\bb\in S^{d-1}(\tau)$. Corollary~\ref{ubcor},
which gives an upper bound on $r(\tau)$,
follows immediately from Proposition~\ref{rrisk} and the definition of
$r(\tau)$.

In practice, the signal strength $\tau= \llVert  \bb\rrVert  $ is
typically unknown. Thus, with $\bb\in S^{d-1}(\tau)$, $\hat{\bb
}_r(\tau)$
may be viewed as an oracle estimator. In cases where the signal
strength is not prespecified,
Proposition~\ref{rrisk} implies that $\hat{\bb}_r(\llVert  \bb\rrVert  )$ is the
oracle estimator with minimal risk among ridge estimators. We will
refer to both $\hat{\bb}_r(\tau)$ and $\hat{\bb}_r(\llVert  \bb\rrVert  )$ as the
oracle ridge estimator, according to whether or not $\bb\in
S^{d-1}(\tau)$ has been specified in advance. In Section~\ref{sec4}, we discuss
adaptive ridge estimators that utilize an estimate of the signal
strength.

Expressions similar to those in Proposition~\ref{rrisk} for the
optimal ridge parameter and the risk ridge estimators have appeared
previously in the literature (see, e.g., the review article by
Draper and Van~Nostrand \cite{draper1979ridge}). However, other
existing results on the risk
of ridge estimators tend
to either (i) be significantly more complex than Proposition~\ref
{rrisk} or (ii) pertain to the Bayes risk of ridge regression,
assuming that $\bb$ follows a normal prior distribution.
Proposition~\ref{rrisk}
is a simple, yet conclusive result for the optimal ridge parameter
with respect to the frequentist risk $R(\hat{\bb},\bb)$. Its
simplicity follows largely from the symmetry in our
formulation of the problem; in particular, we are focusing
on unconditional risk and the distribution of $X$ is
orthogonally invariant.

\subsection{Asymptotic risk}\label{sec2.2}

It appears that the risk formula (\ref{rrisk1}) cannot be further
simplified with ease. However, results
from random matrix theory yield a closed-form expression for the
asymptotic risk. For $\rho\in(0,\infty)$, the Mar\v{c}enko--Pastur
density $f_{\rho}$
is defined by
\[
f_{\rho}(z) = \max \bigl\{ \bigl(1 - \rho^{-1} \bigr), 0 \bigr
\} \d_0(z) + \frac{1}{2\uppi\rho z} \sqrt{(b - z) (z - a)}
\1_{(a,b)}(z),
\]
where $a = (1 - \sqrt{\rho})^2$, $b = (1 + \sqrt{\rho})^2$,
$\d_0(\cdot)$ is the Dirac delta, and $\1_{(a,b)}(\cdot)$ is the indicator
function of the open interval $(a,b)$.
The density $f_{\rho}$ determines the Mar\v{c}enko--Pastur
distribution, which is the limiting distribution of the
eigenvalues of $n^{-1}X^TX$, if $n \to\infty$ and $d/n \to\rho\in
(0,\infty)$ (Mar{\v{c}}enko and Pastur \cite
{marcenko1967distribution}); it also determines the
corresponding cumulative distribution function, $F_{\rho}(t) = \int_{-\infty}^t f_{\rho}(z) \,\mathrm{d}z$. The Stieltjes transform
of the Mar\v{c}enko--Pastur distribution is defined by
%
%
\begin{eqnarray}\label{stieltjes}
m_{\rho}(s) & =& \int\frac{1}{z - s} f_{\rho}(z) \,\mathrm{d}z =
\int\frac
{1}{z - s} \,\mathrm{d}F_{\rho}(z)\nonumber
\\[-8pt]\\[-8pt]
 & =& -\frac{1}{2\rho
s} \bigl\{s + \rho- 1 + \sqrt{(s + \rho-
1)^2 - 4\rho s} \bigr\},\qquad  s < 0.\nonumber
\end{eqnarray}
The main result of this section implies that if
$\bb\in S^{d-1}(\tau)$, then the risk of the oracle ridge estimator
may be approximated by
$(d/n)m_{d/n}\{-d/(n\tau^2)\}$.

%
\begin{thm}\label{mprisk}
Suppose that $0 < \rho_- \leq d/n \leq
\rho_+ < \infty$ for some fixed constants $\rho_-, \rho_+ \in\R$.
\begin{enumerate}[(b)]
\item[(a)] If $0 < \rho_- < \rho_+ < 1$ or $1 < \rho_- <
\rho_+ < \infty$ and $|n-d| > 5$, then
\[
\sup_{\bb\in S^{d-1}(\tau)} \biggl\llvert R \bigl\{\hat{\bb}_r(
\tau),\bb \bigr\} - \frac{d}{n}m_{d/n} \biggl(-\frac{d}{n\tau^2}
\biggr) \biggr\rrvert= \mathrm{O} \biggl(\frac{\tau^2}{1 + \tau^2}n^{-1/2}
\biggr).
\]

\item[(b)] If $0 < \rho_- < 1 < \rho_+ < \infty$, then
\[
\sup_{\bb\in S^{d-1}(\tau)} \biggl\llvert R \bigl\{\hat{\bb}_r(
\tau),\bb \bigr\} - \frac{d}{n}m_{d/n} \biggl(-\frac{d}{n\tau^2}
\biggr) \biggr\rrvert= \mathrm{O} \bigl(\tau^2n^{-1/8} \bigr).
\]
\end{enumerate}
\end{thm}

Theorem~\ref{mprisk} is proved in Appendix \hyperref[appA]{A}. Since
$R\{\hat{\bb}_r(\tau),\bb\}$ is constant over $\bb\in S^{d-1}(\tau
)$, the
supremums in parts (a) and (b) of Theorem~\ref{mprisk} are
somewhat superfluous; however, they serve to emphasize that the upper
bounds do not depend on any particular value of $\bb\in
S^{d-1}(\tau)$.

Let $0 \leq s_d \leq s_{d-1} \leq\cdots\leq s_1$ denote the ordered
eigenvalues of $n^{-1}X^TX$ and define the empirical cumulative
distribution function $\F_{n,d}(s) = d^{-1}\sum_{j = 1}^d
\1_{(-\infty,s_j]}(s)$. There are two keys to the proof of
Theorem~\ref{mprisk}. The first is
the observation that if $\bb\in S^{d-1}(\tau)$, then, by
Proposition~\ref{rrisk},
\[
\frac{n}{d}R \bigl\{\hat{\bb}_r(\tau),\bb \bigr\} =
\frac{1}{d}E \biggl[\tr \biggl\{ \biggl(\frac{1}{n}X^TX +
\frac{d}{n\tau^2}I_d \biggr)^{-1} \biggr\} \biggr] = E
\biggl\{\int\frac{1}{s + d/(n\tau^2)} \,\mathrm{d}\F_{n,d}(s) \biggr\};
\]
in other words, the risk of the oracle ridge estimator
is the expected value of the Stieltjes transform of $\F_{n,d}$. The
second key is Theorem~1.1 of Bai \textit{et al.} \cite
{bai2003convergence}, which states
that under the conditions of Theorem~\ref{mprisk},
%
%
\begin{equation}
\label{bai} \hspace*{-10pt}\sup_{s \in\R} \bigl\llvert E \bigl\{
\F_{n,d}(s) \bigr\} - F_{d/n}(s) \bigr\rrvert= \lleft\{
\begin{array} {l@{\qquad}l} \mathrm{O} \bigl(n^{-1/2} \bigr) &
\mbox{if } %
 0 < \rho_- < \rho_+ < 1
\mbox{ or } 1 < \rho_- < \rho_+ < \infty,  %
\\\noalign{\vspace*{3pt}}
\mathrm{O} \bigl(n^{-1/8} \bigr) & \mbox{if } 0 < \rho_- < 1 < \rho_+ <
\infty. \end{array} %
\rright.
\end{equation}
The different rates in (\ref{bai}) depending on whether or not $\rho_-
< 1 < \rho_+$ helps to explain why these situations are considered
separately in Theorem~\ref{mprisk} above; more fundamentally, the
major difference between the two cases is that if $d/n \to1$
(corresponding to the setting where $\rho_- < 1 < \rho_+$),
then 0 is contained in the support of the continuous part of the Mar\v
{c}enko--Pastur
distribution, which complicates the analysis.

The asymptotic risk of the oracle ridge estimator, when
$d/n \to\rho\in(0,\infty)$, is given
explicitly in the following corollary, which follows immediately from
Theorem~\ref{mprisk}.

%
\begin{cor}\label{corrisk}
For $\rho\in(0,\infty)$ and $\tau\in[0,\infty)$ define the
asymptotic risk of the oracle ridge estimator
\[
R_r(\tau,\rho) = \frac{1}{2\rho} \bigl[\tau^2(\rho-
1) - \rho+ \sqrt{ \bigl\{\tau^2(\rho- 1) - \rho \bigr\}^2
+ 4\rho^2\tau^2} \bigr].
\]
\begin{enumerate}[(b)]
\item[(a)] If $\rho\in(0,\infty)\setminus\{1\}$, then
\[
\lim_{d/n \to\rho} \sup_{\bb\in\R^d} \bigl\llvert R \bigl
\{\hat{\bb}_r\bigl(\llVert \bb\rrVert \bigr),\bb \bigr\} -
R_r\bigl(\llVert \bb\rrVert ,d/n\bigr) \bigr\rrvert= 0.
\]

\item[(b)] If $0 \leq T < \infty$ is a fixed real number, then
\[
\lim_{d/n \to1} \mathop{\sup_{\bb\in\R^d;}}_{ 0 \leq\llVert  \bb\rrVert
\leq T}
\bigl\llvert R \bigl\{\hat{\bb}_r\bigl(\llVert \bb\rrVert \bigr),\bb
\bigr\} - R_r\bigl(\llVert \bb\rrVert ,d/n\bigr) \bigr\rrvert= 0.
\]
\end{enumerate}
\end{cor}

In Corollary~\ref{corrisk} and throughout the paper, the notation
$\lim_{d/n
\to\rho}$ indicates the limit as $n \to\infty$ and $d/n \to
\rho$. Corollary~\ref{corrisk} implies that if $d/n \to\rho\in
(0,\infty)\setminus\{1\}$, then
the risk of the oracle ridge estimator $\hat{\bb}_r(\llVert  \bb\rrVert  )$
converges to the
asymptotic risk $R_r(\llVert  \bb\rrVert  ,d/n)$ uniformly over all $\bb\in\R^d$;
if $d/n \to1$, then the
convergence is uniform over compact sets.

It is clear from Theorem~\ref{mprisk} and Corollary~\ref{corrisk} that
if $d/n \to\rho\in(0,\infty)$, then the spectral distribution of
$n^{-1}X^TX$ plays a prominent role in determining the risk of the
oracle ridge estimator via the Mar\v{c}enko--Pastur law; if $d/n \to0$
or $d/n \to\infty$, then
its role subsides, as illustrated by the following proposition.

%
\begin{prop}\label{r0inf}
\begin{enumerate}[(b)]
\item[(a)] [$d/n \to0$] For $\rho, \tau\in[0,\infty)$ define
%
%
\begin{equation}
\label{r0inf1} R_r^0(\tau,\rho) = \frac{\rho\tau^2}{\rho+ \tau^2}.
\end{equation}
Then
\[
\lim_{d/n \to0} \sup_{\bb\in\R^d} \biggl\llvert
\frac{R\{\hat{\bb
}_r(\llVert  \bb\rrVert  ),\bb\} }{R_r^0(\llVert  \bb\rrVert  ,d/n)} - 1 \biggr\rrvert= 0.
\]

\item[(b)] [$d/n\to\infty$] Let $0 \leq T < \infty$
be a fixed real number. Then
\[
\lim_{d/n \to\infty} \mathop{\sup_{\bb\in\R^d;}}_{ 0 \leq
\llVert  \bb\rrVert   \leq T}
\bigl\llvert R \bigl\{\hat{\bb}_r\bigl(\llVert \bb\rrVert \bigr),\bb
\bigr\} - \llVert \bb\rrVert ^2 \bigr\rrvert= 0.
\]
\end{enumerate}
\end{prop}

Proposition~\ref{r0inf} is proved in Appendix \hyperref[appA]{A}. It gives the
asymptotic risk of the oracle ridge
estimator in settings where $d/n \to0$ and $\infty$. Expressions
like (\ref{r0inf1}) are common in the analysis of linear estimators
for the Gaussian sequence model (Pinsker \cite
{pinsker1980optimal}). Thus, if
$d/n \to0$, then features of $R\{\hat{\bb}_r(\llVert  \bb\rrVert  ),\bb\} $
deriving from the
random predictors $X$ are less apparent.

Now consider the null
estimator $\hat{\bb}_{\mathrm{null}} = 0$ and notice that $R(\hat{\bb
}_{\mathrm{null}},\bb) =
\llVert  \bb\rrVert  ^2$. Proposition~\ref{r0inf}(b) implies that
if $d/n \to\infty$, then the oracle ridge estimator is asymptotically
equivalent to $\hat{\bb}_{\mathrm{null}}$. In Section~\ref{sec3}, we argue that if
$d/n \to\infty$, then $\hat{\bb}_{\mathrm{null}}$ is in fact asymptotically
minimax for the problem (\ref{mr}). In other words, non-trivial
estimation is impossible in (\ref{mr})
when $d/n
\to\infty$.

Combined with Theorem~\ref{mprisk}, Proposition~\ref{r0inf} implies
that the asymptotic
risk of the oracle ridge estimator $R_r(\tau,\rho)$ extends
continuously to $\rho= 0$ and $\rho= \infty$. For $\tau\geq0$, we
define $R_r(\tau,0)
= 0$ and $R_r(\tau,\infty) = \tau^2$.

\section{An equivalent Bayes problem}\label{sec3}

In this section, we use an equivariance argument to reduce the minimax
problem (\ref{mr}) to an equivalent Bayes problem. We then show that
ridge regression solves the Bayes problem, asymptotically.

\subsection{The uniform measure on \texorpdfstring{$S^{d-1}(\tau)$}{$S^{d-1}(tau)$} and equivariance}\label{sec3.1}

Let $\pi_{S^{d-1}(\tau)}$ denote the uniform measure on $S^{d-1}(\tau
)$. Define
the Bayes risk
%
%
\begin{equation}
\label{bayes} r_B(\tau) = \inf_{\hat{\bb}} \int
_{S^{d-1}(\tau)} R(\hat{\bb},\bb) \,\mathrm{d}\pi_{S^{d-1}(\tau)}(\bb) = \inf
_{\hat{\bb}} E_{\pi
_{S^{d-1}(\tau)}} \bigl(\llVert \hat{\bb} - \bb\rrVert
^2 \bigr),
\end{equation}
where the expectation $ E_{\pi_{S^{d-1}(\tau)}}$ is taken with
respect to the joint
distribution of $(X,\ee,\bb)$, with $\bb\sim\pi_{S^{d-1}(\tau)}$
independent of
$(X,\ee)$. The Bayes estimator
\[
\hat{\bb}_{S^{d-1}(\tau)} = E_{\pi_{S^{d-1}(\tau)}}(\bb|\y,X)
\]
satisfies
%
%
\begin{equation}
\label{bayesSoln} r_B(\tau) = E_{\pi_{S^{d-1}(\tau)}} \bigl\{ \llVert \hat{
\bb }_{S^{d-1}(\tau)} - \bb\rrVert ^2 \bigr\}.
\end{equation}

Let $\O(d)$ denote the group of $d \times d$ orthogonal matrices. As
with $\ee$ and $X$, the distribution
$\pi_{S^{d-1}(\tau)}$ is invariant under orthogonal
transformations; that is, if $U \in\O(d)$ and $\bb\sim
\pi_{S^{d-1}(\tau)}$, then $U\bb\sim\pi_{S^{d-1}(\tau)}$. A corresponding
feature of the estimator $\hat{\bb}_{S^{d-1}(\tau)}$ is that it is
\emph{equivariant} with respect to orthogonal transformations.

%
\begin{defn}\label{ee}
An estimator $\hat{\bb} = \hat{\bb}(\y,X)$ is \emph{orthogonally
equivariant} if
\[
\hat{\bb}(\y,XU) = U^T\hat{\bb}(\y,X)
\]
for all $d \times d$ orthogonal matrices $U \in\mathcal{O}(d)$.
\end{defn}

Let
\[
\E= \E_{d,n} = \{\hat{\bb}; \hat{\bb} \mbox{ is an orthogonally
equivariant estimator for } \bb\}.
\]
Then one easily checks that $\hat{\bb}_{S^{d-1}(\tau)} \in\E$.
Additionally, notice that $\hat{\bb}_r(\tau) \in\E$ is
orthogonally equivariant. The following proposition is proved in
Appendix \hyperref[appA]{A}.

%
\begin{prop}\label{mrequiv}Suppose that $\tau\geq0$ and that
$\bb_1,\bb_2 \in S^{d-1}(\tau)$.
\begin{enumerate}[(b)]
\item[(a)] If $\hat{\bb}$ is an orthogonally
equivariant estimator, then the risk of $\hat{\bb}$ is constant over
$S^{d-1}(\tau)$; that is, $R(\hat{\bb},\bb_1) =
R(\hat{\bb},\bb_2)$.

\item[(b)]
\[
\hspace*{-10pt}r(\tau) = \inf_{\hat{\bb} \in\E} \sup_{\bb\in S^{d-1}(\tau)} R(\hat{
\bb},\bb) = R\{\hat{\bb}_{S^{d-1}(\tau)},\bb_1\} = E_{\pi
_{S^{d-1}(\tau)}}
\bigl\{\llVert \hat{\bb}_{S^{d-1}(\tau)} - \bb\rrVert ^2 \bigr\} =
r_B(\tau).
\]
\end{enumerate}
\end{prop}

Proposition~\ref{mrequiv}(a) implies that all orthogonally
equivariant estimators
have constant risk over spheres $S^{d-1}(\tau)$; we first noted that ridge
regression possesses this property in a remark following
Proposition~\ref{rrisk}. Proposition~\ref{mrequiv}(b) implies that
the Bayes
problem (\ref{bayes})
and the minimax problem (\ref{mr}) are equivalent. Proposition~\ref
{mrequiv}(b) also implies that the estimator
$\hat{\bb}_{S^{d-1}(\tau)}$ is minimax over $S^{d-1}(\tau)$. While
this, in a sense,
``solves'' the main problem of interest (\ref{mr}), there are several
caveats. For instance, the estimator
$\hat{\bb}_{S^{d-1}(\tau)}$ is an oracle estimator (it depends on
$\tau$)
and is difficult to compute, even if $\tau$ is known. Furthermore,
Proposition~\ref{mrequiv} provides no information about the magnitude of
$r(\tau)$. In the next section,
we show that when $d$ is large, $r(\tau) =
R\{\hat{\bb}_{S^{d-1}(\tau)},\bb\} \approx R\{\hat{\bb}_r(\tau
),\bb\}
\approx R_r(\tau,\rho)$ for $\bb\in S^{d-1}(\tau)$. In addition to
providing quantitative information about $r(\tau)$, this result
suggests that
ridge regression may be an appealing alternative to
$\hat{\bb}_{S^{d-1}(\tau)}$, especially when combined with results on
adaptive ridge estimators in Section~\ref{sec4}.

\subsection{Ridge regression and asymptotic optimality}\label{sec3.2}

Recall that the minimax estimator $\hat{\bb}_{S^{d-1}(\tau)}$ is the
posterior mean of $\bb$, under the assumption that $\bb\sim\pi
_{S^{d-1}(\tau)}$ is
uniformly distributed on the sphere $S^{d-1}(\tau)$. On the other
hand, the
oracle ridge estimator $\hat{\bb}_r(\tau) = E_{N(0,\tau
^2/dI_d)}(\bb|\y,X)$
may be interpreted
as the posterior mean of $\bb$
under the assumption that $\bb\sim N(0,\tau^2/dI_d)$ is normally
distributed and independent of $(X,\ee)$. If $d$ is
large, then the normal distribution $N(0,\tau^2/dI_d)$ is ``close'' to
the uniform distribution on $S^{d-1}(\tau)$ (there is an enormous body
of literature that makes this idea
more precise -- Diaconis and Freedman \cite
{diaconis1987dozen} attribute early work to
Borel \cite{borel1914introduction} and
L{\'e}vy \cite{levy1922lecons}). Thus, it is
reasonable to expect that if $d$ is large and $\bb\in S^{d-1}(\tau)$,
then $\hat{\bb}_{S^{d-1}(\tau)} \approx\hat{\bb}_r(\tau)$
and that the two estimators have similar risk properties. This is
the content of the main
result in this section, which is
essentially a multivariate extension of Theorem~3.1 from Marchand \cite{marchand1993estimation}.

%
\begin{thm} \label{main} Suppose that $n > 2$ and let $s_1 \geq\cdots
\geq s_{d \wedge n} > 0$ denote the nonzero (with probability $1$)
eigenvalues of $n^{-1}X^TX$. Let $\tau\geq0$.
\begin{enumerate}[(b)]
\item[(a)] If $d \leq n$ and $\bb\in S^{d-1}(\tau)$, then
\[
R\{\hat{\bb}_{S^{d-1}(\tau)},\bb\} \leq R \bigl\{\hat{\bb}_r(\tau),
\bb \bigr\} \leq R\{\hat{\bb}_{S^{d-1}(\tau)},\bb\} + \frac
{1}{d}E \biggl[
\frac{s_1}{s_d}\tr \biggl\{ \biggl(X^TX + \frac{d}{\tau^2}I_d
\biggr)^{-1} \biggr\} \biggr].
\]

\item[(b)] If $d > n$ and $\bb\in S^{d-1}(\tau)$, then
\begin{eqnarray*}
\hspace*{-10pt}R\{\hat{\bb}_{S^{d-1}(\tau)},\bb\} \leq R \bigl\{\hat{\bb}_r(\tau),
\bb \bigr\} & \leq& R\{\hat{\bb}_{S^{d-1}(\tau)},\bb\} + \frac
{1}{n}E \biggl[
\frac{s_1}{s_n}\tr \biggl\{ \biggl(XX^T + \frac{d}{\tau^2}I_n
\biggr)^{-1} \biggr\} \biggr]
\\
&&{} + \frac{2(d-n)}{\tau^2(n-2)}E \biggl[\tr \biggl\{ \biggl(XX^T +
\frac{d}{\tau^2}I_n \biggr)^{-2} \biggr\} \biggr].
\end{eqnarray*}
\end{enumerate}
\end{thm}

Theorem~\ref{main} is proved in Appendix \hyperref[appB]{B}. The bound
$R\{\hat{\bb}_{S^{d-1}(\tau)},\bb\} \leq R\{\hat{\bb}_r(\tau
),\bb\}$ follows
immediately from Proposition~\ref{mrequiv}(b) and Corollary~\ref
{ubcor}. Proving the
required upper bounds on
$R\{\hat{\bb}_r(\tau),\bb\}$ (which, by Proposition~\ref{mrequiv}(b), are equivalent
to lower bounds on $r(\tau)$) is fairly complex and
involves transforming the linear model into an equivalent
sequence
model, along with the application of classical information identities
(Brown \cite{brown1971admissible}) and
inequalities (Stam \cite{stam1959some}). In
the remainder of this section, we discuss some of the implications of
Theorem~\ref{main}.

Asymptotically, Theorem~\ref{main} is primarily significant for
settings where $d/n \to\rho\in(0,\infty)$. If $d/n \to\rho\in
(0,\infty)\setminus\{1\}$, then the upper bounds in Theorem~\ref{main}
are $\mathrm{O}(n^{-1})$ and $R\{\hat{\bb}_r(\tau),\bb\} \approx
R\{\hat{\bb}_{S^{d-1}(\tau)},\bb\} = r(\tau)$, where $\bb\in
S^{d-1}(\tau)$;
by Corollary~\ref{corrisk}, we can further conclude that $r(\tau)
\approx R_r(\tau,d/n)$. The
case where $d/n \to1$ is somewhat problematic, because then
$E(s_d^{-1}) \to\infty$; however, some conclusions can be made in
this case by continuity arguments, for example, Corollary~\ref
{rhoposcor}(b) below.

%
\begin{prop} \label{rhopos}
Suppose that $0 < \rho_- \leq d/n \leq
\rho_+ < \infty$ for some fixed constants $\rho_-, \rho_+ \in\R$ and
that $0 < \rho_- < \rho_+ < 1$ or $1 < \rho_- <
\rho_+ < \infty$. If $|n-d| > 5$, then
\begin{eqnarray*}
\sup_{\bb\in S^{d-1}(\tau)} \bigl\llvert R \bigl\{\hat{\bb}_r(
\tau),\bb \bigr\} - R\{\hat{\bb}_{S^{d-1}(\tau)},\bb\} \bigr\rrvert& =& \sup
_{\bb\in
S^{d-1}(\tau)} \bigl\llvert R \bigl\{\hat{\bb}_r(\tau),\bb
\bigr\} - r(\tau) \bigr\rrvert
\\
& =& \mathrm{O} \biggl(\frac{\tau^2}{\tau^2
+ 1}n^{-1} \biggr).
\end{eqnarray*}
\end{prop}

%
\begin{cor} \label{rhoposcor} Let $R_r(\tau,\rho)$ be the asymptotic
risk of the ridge estimator defined in Corollary~\ref{corrisk}.
\begin{enumerate}[(b)]
\item[(a)] If
$\rho\in(0,\infty)\setminus\{1\}$, then
\[
\lim_{d/n\to\rho} \sup_{0 \leq\tau} \bigl\llvert
R_r(\tau,d/n) - r(\tau) \bigr\rrvert= \lim_{d/n\to\rho}
\sup_{\bb\in\R^d} \bigl\llvert R \bigl\{\hat{\bb}_r\bigl(
\llVert \bb\rrVert \bigr),\bb \bigr\} - r\bigl(\llVert \bb\rrVert \bigr) \bigr
\rrvert= 0.
\]

\item[(b)] If $0 \leq T < \infty$ is a fixed real number, then
\[
\lim_{d/n \to1} \sup_{0 \leq\tau\leq T} \bigl\llvert
R_r(\tau,d/n) - r(\tau) \bigr\rrvert= \lim_{d/n \to1}
\mathop{\sup_{\bb\in\R^d;}}_{ 0
\leq\llVert  \bb\rrVert   \leq T} \bigl\llvert R \bigl\{\hat{
\bb}_r\bigl(\llVert \bb\rrVert \bigr),\bb \bigr\} - r\bigl(\llVert
\bb\rrVert \bigr) \bigr\rrvert= 0.
\]
\end{enumerate}
\end{cor}

Proposition~\ref{rhopos} follows directly from Theorem~\ref{main} and
Lemma~\ref{c2} (found in Appendix \hyperref[appC]{C}). Corollary~\ref{rhoposcor}(a)
follows immediately from Proposition~\ref{rhopos} and
Corollary~\ref{corrisk}(a). Corollary~\ref{rhoposcor}(b)
may be proved similarly to part (a), while making use of the inequality
$r_{d,n-k}(\tau) \leq r_{d,n}(\tau)$ for integers $0 \leq k < n$ in
order to avoid issues around $d/n \approx1$. Corollary~\ref
{rhoposcor} implies that if $d/n \to\rho\in
(0,\infty)$, then the minimax risk
$r(\tau)$ is asymptotically equivalent to the asymptotic risk of the
oracle ridge estimator and that the oracle ridge estimator is
asymptotically minimax.

Corollary~\ref{rhoposcor} also
provides the means for relating the minimax problem over
$\ell^2$-spheres (\ref{mr}) to the minimax problem over $\ell^2$-balls
(\ref{mrb}). Since $S^{d-1}(\tau) \subseteq B_d(\tau)$, we have
$r(\tau) \leq
\bar{r}(\tau)$. Furthermore, one easily checks that
\[
\sup_{\bb\in B_d(\tau)}R \bigl\{\hat{\bb}_r(\tau),\bb \bigr\}
= \sup_{\bb
\in
S^{d-1}(\tau)} R \bigl\{\hat{\bb}_r(\tau),\bb \bigr
\}.
\]
Thus, if $d/n \to\rho\in
(0,\infty)$, then
%
%
\begin{equation}
\label{mrbmr} r(\tau) \leq\bar{r}(\tau) \leq\sup_{\bb\in
B_d(\tau)} R
\bigl\{\hat{\bb}_r(\tau),\bb \bigr\} = \sup_{\bb\in
S^{d-1}(\tau)}R
\bigl\{ \hat{\bb}_r(\tau),\bb \bigr\} \to r(\tau).
\end{equation}
It follows that if $d/n \to\rho\in(0,\infty)$, then the minimax
risk over $S^{d-1}(\tau)$ is equivalent to the minimax risk over
$B_d(\tau)$
and that the ridge estimator $\hat{\bb}_r(\tau)$
is asymptotically minimax for both problems.

When $d/n \to0$
or $\infty$, asymptotics for the minimax risk $r(\tau)$ are more
straightforward. The following proposition summarizes the behavior of
$r(\tau)$ in these settings.

%
\begin{prop} \label{rhoalt}
\begin{enumerate}[(b)]
\item[(a)] [$d/n \to0$] Let
$R_r^0(\tau,\rho)$ be the risk function (\ref{r0inf1}). Then
\[
\mathop{\lim_{d/n \to0}}_{ d\to\infty} \sup_{\tau\geq0}
\biggl\llvert\frac{R_r^0(\tau,d/n)}{r(\tau)} - 1 \biggr\rrvert= 0.
\]

\item[(b)] [$d/n \to\infty$] Let $0 < T< \infty$ be fixed. Then
\[
\lim_{d/n \to\infty} \sup_{0 \leq\tau\leq T} \bigl\llvert r(\tau)
- \tau^2 \bigr\rrvert= 0.
\]
\end{enumerate}
\end{prop}

Proposition~\ref{rhoalt}(a) is a straightforward consequence of
Theorem~\ref{main}, Proposition~\ref{r0inf}, and Lem\-ma~\ref{c2}.
Proposition~\ref{rhoalt}(b) follows from general properties of
orthogonally equivariant estimators; in particular, one can check that
if $d \geq n$, then
\[
R(\hat{\bb},\bb) \geq\frac{d-n}{d}\llVert \bb\rrVert ^2
\]
for all orthogonally equivariant estimators $\hat{\bb}$.

Proposition~\ref{rhoalt} gives precise
asymptotics for $r(\tau)$ when $d \to\infty$ and $d/n \to0$ or
$\infty$. While Proposition~\ref{rhoalt} does not directly reference
the ridge estimator, combined
with Proposition~\ref{r0inf} it implies that $\hat{\bb}_r(\tau)$
is asymptotically optimal for the minimax problem (\ref{mr}) when
$d \to
\infty$ and $d/n
\to0$ or $\infty$. Note that the
null estimator $\hat{\bb}_{\mathrm{null}} = 0$ is also asymptotically optimal
for (\ref{mr}) when $d/n \to\infty$. We point out that the
condition $d \to\infty$ in Proposition~\ref{rhoalt}(a) appears to be
necessary, as it drives the approximation $\pi_{S^{d-1}(\tau)}
\approx
N(0,\tau^2/dI_d)$ underlying Theorem~\ref{main}.

\section{An adaptive ridge estimator}\label{sec4}

To this point, we have focused on the oracle ridge estimator
$\hat{\bb}_r(\tau)$, where $\tau= \llVert  \bb\rrVert  $ is the signal strength.
Typically, $\tau$ is unknown and, consequently,
$\hat{\bb}_r(\tau)$ is non-implementable. A natural strategy is to
replace $\tau$ with an estimate, $\hat{\tau}$.

Define
%
%
\begin{equation}
\label{tauhat} \hat{\tau}^2 = \max \biggl\{\frac{1}{n}\llVert
\y\rrVert ^2 - 1,0 \biggr\}
\end{equation}
and define the adaptive
ridge estimator
\[
\check{\bb}_r = \hat{\bb}_r(\hat{\tau}).
\]
Observe that $\check{\bb}_r \in\E$ is orthogonally equivariant. One
can check that $\sup_{\bb\in S^{d-1}(\tau)}
E_{\bb}(\hat{\tau}/\tau- 1)^2 \to0$ whenever $n \to\infty$ (see
Lemma~\ref{c5}); thus, $\hat{\tau}$ is a
reasonable estimator for $\tau$. The next result relates the
risk of the adaptive
ridge estimator $\check{\bb}_r$ to that of the oracle ridge
estimator. It is proved in Appendix \hyperref[appA]{A}.

%
\begin{thm}\label{adapt}
Suppose that $0 < \rho_- <d/n < \rho_+ < \infty$, where
$\rho_-,\rho_+\in\R$ are fixed constants satisfying $0 < \rho_- <
\rho_+ < 1$ or $1 < \rho_- < \rho_+ <
\infty$. Also suppose that $|n -
d| > 9$ and $n > 8$. Then
%
%
\begin{equation}
\label{adapt1} \sup_{\bb\in S^{d-1}(\tau)} \bigl\llvert R(\check{
\bb}_r,\bb)-R \bigl\{\hat{\bb}_r(\tau),\bb \bigr\} \bigr
\rrvert= \mathrm{O} \biggl(\frac{1}{\tau^2 + 1}n^{-1/2} \biggr)
\end{equation}
and
%
%
\begin{equation}
\label{adapt2} \sup_{\bb\in\R^d} \bigl\llvert R(\check{
\bb}_r,\bb) - R_r\bigl(\llVert \bb\rrVert ,d/n\bigr)
\bigr\rrvert= \mathrm{O} \bigl(n^{-1/2} \bigr),
\end{equation}
where $R_r(\tau,\rho)$ is the asymptotic risk of the oracle ridge
estimator defined in Corollary~\ref{corrisk}.
\end{thm}

Theorem~\ref{adapt} implies that if $d/n \to\rho\in
(0,\infty)\setminus\{1\}$, then the risk of the adaptive ridge
estimator converges uniformly to that of the oracle ridge estimator
and its asymptotic risk is given explicitly by $R_r(\tau,\rho)$.
If $d/n \to\rho\in
(0,\infty)\setminus\{1\}$ and $\tau^2 = \llVert  \bb\rrVert  ^2 \gg n^{-1/2}$,
then it
follows from Theorem~\ref{adapt} that $R(\check{\bb}_r,\bb)/R\{\hat
{\bb}_r(\tau),\bb
\}
\to1$. On the other hand, if $d/n \to\rho\in
(0,\infty)\setminus\{1\}$ and $\tau^2 = \mathrm{O}(n^{-1/2})$, then
$R\{\hat{\bb}_r(\tau),\bb\} = \mathrm{O}(n^{-1/2})$ and the
limit of $R(\check{\bb}_r,\bb)/R\{\hat{\bb}_r(\tau),\bb\}$ does
not follow
readily from Theorem~\ref{adapt}. In other words, the effectiveness
of the adaptive ridge estimator is less clear when $\tau^2 = \llVert  \bb
\rrVert  ^2$ is very
small.

If $d/n \to0$ or $d/n \to1$, then results similar to Theorem~\ref
{adapt} may be
obtained for the adaptive ridge estimator, but the results are more
delicate; results for $d/n \to
\infty$ are, in a sense, unnecessary because the oracle ridge
estimator is equivalent to $\hat{\bb}_{\mathrm{null}}$ in this setting. If
$d/n \to0$, then the relevant quantity is the risk ratio
$R(\check{\bb}_r,\bb)/R\{\hat{\bb}_r(\tau),\bb\}$, rather than
the risk
difference considered in Theorem~\ref{adapt}, and one must carefully
track the magnitude of $\tau^2 = \llVert  \bb\rrVert  ^2$ relative to $d/n$. Ultimately,
however, when $d/n \to0$ the message is the same as the case where
$d/n \to\rho\in(0,\infty)\setminus\{1\}$: If $\tau^2$ is not
too small, then the adaptive ridge estimator performs nearly as well
as the oracle ridge estimator. If $d/n \to1$, then the rate in
(\ref{adapt1}) may be different, depending on
$\tau^2$ and the magnitude of $|d/n - 1|$, for example,
Bai \textit{et al.} \cite{bai2003convergence}.

\subsection{Adaptive minimax estimation}\label{sec4.1}

Theorem~\ref{adapt} compares the risk of the adaptive ridge estimator to
that of the oracle ridge estimator. The next result, which follows
immediately from Theorem~\ref{adapt} and Proposition~\ref{rhopos},
compares the risk of the adaptive ridge estimator to $r(\tau)$.

%
\begin{prop} \label{adaptmr}
Suppose that $\rho_-,\rho_+ \in\R$ are fixed constants satisfying $0
< \rho_- < \rho_+ < 1$ or $1 < \rho_- < \rho_+ < \infty$. Suppose
further that $0 < \rho_- \leq d/n \leq
\rho_+ < \infty$. If $|n-d| > 9$ and $n > 8$, then
%
%
\begin{equation}
\label{adaptmr1} \sup_{\bb\in S^{d-1}(\tau)} \bigl\llvert R(\check{
\bb}_r,\bb) - r(\tau) \bigr\rrvert= \mathrm{O} \biggl(
\frac{\tau^2}{\tau^2
+ 1}n^{-1} \biggr) + \mathrm{O} \biggl(\frac{1}{\tau^2 + 1}
n^{-1/2} \biggr).
\end{equation}
\end{prop}

Combined with Proposition~\ref{rhopos}, Proposition~\ref{adaptmr}
implies that if $d/n \to\rho\in
(0,\infty)\setminus\{1\}$, then $\check{\bb}_r$ is adaptive asymptotic
minimax over spheres $S^{d-1}(\tau)$, provided $\tau^2 \gg n^{-1/2}$.

\subsection{Equivariance}\label{sec4.2}

In Section~\ref{sec3.1}, we discussed connections between the minimax problem
(\ref{mr}) and equivariance. Previously in this section, we noted that
the adaptive ridge estimator $\check{\bb}_r$ is orthogonally
equivariant and adaptive asymptotic minimax over spheres $S^{d-1}(\tau
)$. The following
is an asymptotic optimality result for $\check{\bb}_r$, which pertains
to the class of orthogonally equivariant estimators $\E$.

%
\begin{prop} \label{Oadapt}
Suppose that $\rho\in(0,\infty)\setminus\{1\}$. Then
\[
\lim_{d/n \to\rho} \sup_{\bb\in\R^d} \Bigl\llvert R(
\check{\bb}_r,\bb) - \inf_{\hat{\bb} \in\E} R(\hat{\bb},\bb)
\Bigr\rrvert= 0.
\]
\end{prop}

By Proposition~\ref{mrequiv}, $\inf_{\hat{\bb} \in\E} R(\hat{\bb
},\bb) =
r(\llVert  \bb\rrVert  )$. Thus, Proposition~\ref{Oadapt} is a direct consequence
of Proposition~\ref{adaptmr}. Proposition~\ref{Oadapt} implies that
if $d/n \to\rho\in(0,\infty)\setminus\{1\}$, then the adaptive
ridge estimator
$\check{\bb}_r$ is asymptotically optimal among all orthogonally
equivariant estimators. Note that the caveats discussed after the
statement of Theorem~\ref{adapt} relating to small $\llVert  \bb\rrVert  $ also
apply to Proposition~\ref{Oadapt}. More
specifically, if $\llVert  \bb\rrVert   = \mathrm{O}(n^{-1/2})$, then the ratio
$R(\check{\bb}_r,\bb)/\{\inf_{\hat{\bb} \in\E} R(\hat{\bb
},\bb)\}$ is
more relevant than the risk difference considered in Proposition~\ref
{Oadapt} and the precise asymptotic behavior of this ratio is
less clear.
\begin{appendix}
\setcounter{subsection}{1}
\renewcommand{\theequation}{\arabic{equation}}
\section*{Appendix A}\label{appA}

This appendix contains proofs of results stated in the main text, with
the exception of Theorem~\ref{main}; a proof of Theorem~\ref{main}
may be found in Appendix \hyperref[appB]{B}.

\begin{pf*}{Proof of Proposition~\ref{rrisk}} Fix $t \in[0,\infty]$ and
suppose that $\bb\in S^{d-1}(\tau)$. Then
%
%
\begin{eqnarray}\label{prop1a}
\nonumber
R \bigl\{\hat{\bb}_r(t),\bb \bigr\} & =& E_{\bb}
\bigl\{\bigl\llVert \hat{\bb}_r(t) - \bb\bigr\rrVert ^2
\bigr\}\nonumber
\\
& =& E \bigl\{  \bigl\llVert d \bigl(t^2X^TX
+ dI_d \bigr)^{-1}\bb- t^2
\bigl(t^2X^TX + dI_d \bigr)^{-1}X^T
\ee  \bigr\rrVert^2 \bigr\}
\\
 & =& E \bigl\{ \bigl\llVert  d \bigl(t^2X^TX
+ d I_d \bigr)^{-1}\bb \bigr\rrVert
^2 \bigr\} + E \bigl\{  \bigl\llVert t^2
\bigl(t^2X^TX + d I_d \bigr)^{-1}X^T
\ee  \bigr\rrVert^2 \bigr\}.\nonumber
\end{eqnarray}
Since $X$ is orthogonally invariant (i.e., $X$ and $XU$ have the same
distribution for any $U \in\O(d)$), it
follows that
\begin{eqnarray*}
E \bigl\{  \bigl\llVert d \bigl(t^2X^TX + d
I_d \bigr)^{-1}\bb  \bigr\rrVert^2
\bigr\} & =& d^2E \bigl\{\bb^T \bigl(t^2X^TX
+ dI_d \bigr)^{-2}\bb \bigr\}
\\
& =& d^2\tau^2E \bigl\{\mathbf{e}_k^T
\bigl(t^2X^TX + dI_d \bigr)^{-2}
\mathbf{e}_k \bigr\},
\end{eqnarray*}
where $\mathbf{e}_k = (0,\ldots,0,1,0,\ldots,0)^T \in\R^d$ is the
$k$th standard basis vector. Summing over $k = 1,\ldots,d$ above and
dividing by $d$, we obtain
%
%
\begin{equation}
\label{prop1b} E \bigl\{  \bigl\llVert d \bigl(t^2X^TX
+ d I_d \bigr)^{-1}\bb  \bigr\rrVert
^2 \bigr\} = d \tau^2 E \bigl[\tr \bigl\{
\bigl(t^2X^TX + d I_d \bigr)^{-2}
\bigr\} \bigr].
\end{equation}
Additionally, it is clear that
\[
E \bigl\{  \bigl\llVert t^2 \bigl(t^2X^TX
+ d I_d \bigr)^{-1}X^T\ee  \bigr
\rrVert^2 \bigr\} = E \bigl[\tr \bigl\{t^4
\bigl(t^2X^TX + d I_d \bigr)^{-2}X^TX
\bigr\} \bigr].
\]
Combining this with (\ref{prop1a}) and (\ref{prop1b}) yields
\begin{eqnarray*}
R \bigl\{\hat{\bb}_r(t),\bb \bigr\} & =& d \tau^2 E
\bigl[ \tr \bigl\{ \bigl(t^2X^TX + d I_d
\bigr)^{-2} \bigr\} \bigr] + E \bigl[\tr \bigl\{t^4
\bigl(t^2X^TX + d I_d \bigr)^{-2}X^TX
\bigr\} \bigr]
\\
& =& E \bigl[\tr \bigl\{ \bigl(t^2X^TX + d
I_d \bigr)^{-2} \bigl(t^4X^TX + d
\tau^2I_d \bigr) \bigr\} \bigr].
\end{eqnarray*}
Now let $s_1 \geq\cdots\geq s_d \geq0$ denote the eigenvalues of
$n^{-1}X^TX$. Then
\begin{eqnarray*}
R \bigl\{\hat{\bb}_r(t),\bb \bigr\} & =& E \Biggl\{\sum
_{j = 1}^d \frac{t^4ns_j
+ d
\tau^2}{ (t^2ns_j + d )^{2}} \Biggr\}
\\
& =& E \Biggl[\sum_{j = 1}^d \biggl\{
\frac{\tau^2}{\tau^2ns_j + d} + \frac{dns_j(\tau^2 -
t^2)^2}{(t^2ns_j + d)^2(\tau^2ns_j + d)} \biggr\} \Biggr].
\end{eqnarray*}
Clearly, the right-hand
side above is minimized by taking $t = \tau$ and
$R\{\hat{\bb}_r(\tau),\bb\} = E[\tr\{(X^TX + d/\tau^2I_d)^{-1}\}
]$.
\end{pf*}

\begin{pf*}{Proof of Theorem~\ref{mprisk}} Suppose that $\bb\in
S^{d-1}(\tau)$ and let $\F_{n,d}$ be the empirical cumulative
distribution function of the eigenvalues of $n^{-1}X^TX$. Using
integration by parts, for $c \geq0$,
%
%
\begin{eqnarray}\label{mpriska}
\frac{n}{d}\tr \bigl\{ \bigl(X^TX + d/\tau^2
I_d \bigr)^{-1} \bigr\} & =& \int_0^{\infty}
\frac{1}{s + d/(n\tau^2)} \,\mathrm{d}\F_{n,d}(s)
\nonumber
\\
& =& \int_0^c \frac{1}{s + d/(n\tau^2)} \,\mathrm{d}\F_{n,d}(s) + \frac{1}{c +
d/(n\tau^2)} \bigl\{1 - \F_{n,d}(c) \bigr\}
\\
&&{} - \int_c^{\infty} \frac{1}{\{s + d/(n\tau^2)\}^2} \bigl\{1 - \F
_{n,d}(s) \bigr\} \,\mathrm{d}s. \nonumber
\end{eqnarray}
Similarly,
\begin{eqnarray}\label{mpriskb}
m_{d/n} \bigl\{-d/ \bigl(n\tau^2 \bigr) \bigr\} &
= &\int_0^c \frac{1}{s + d/(n\tau^2)}
\,\mathrm{d}F_{d/n}(s) + \frac{1}{c + d/(n\tau^2)} \bigl\{1 - F_{d/n}(c) \bigr\}
\nonumber
\\[-8pt]\\[-8pt]
&&{} - \int_c^{\infty} \frac{1}{\{s + d/(n\tau^2)\}^2} \bigl\{1 -
F_{d/n}(s) \bigr\} \,\mathrm{d}s. \nonumber
\end{eqnarray}

Now let $\Delta= | R\{\hat{\bb}_r(\tau),\bb\} -
(d/n)m_{d/n}\{-d/(n\tau^2)\}|$. Taking $c = 0$ in (\ref{mpriska}) and
(\ref{mpriskb}) implies
\begin{eqnarray*}
\Delta& \leq&\frac{d}{n}\int_0^{\infty}
\frac{1}{\{s +
d/(n\tau^2)\}^2} \bigl\llvert E \bigl\{\F_{n,d}(s) \bigr\} -
F_{d/n}(s) \bigr\rrvert \,\mathrm{d}s
\\
& \leq&\llVert \bb\rrVert ^2 \sup_{s \geq0} \bigl
\llvert E \bigl\{ \F_{n,d}(s) \bigr\} - F_{d/n}(s) \bigr\rrvert,
\end{eqnarray*}
where we have used the fact that $\F_{n,d}(0) = F_{d/n}(0) =
\max\{1-n/d,0\}$, with probability 1.
Thus, it follows from Theorem~1.1 of Bai \textit{et al.}
\cite{bai2003convergence} (see
equation (\ref{bai}) in Section~\ref{sec2.2} above) that
\[
\Delta= \lleft\{ %
\begin{array} {l@{\qquad}l} \mathrm{O} \bigl(
\tau^2 n^{-1/2} \bigr) & \mbox{if } 0 < \rho_- < \rho_+ < 1
\mbox{ or } 1 < \rho_- < \rho_+ < \infty,
\\\noalign{\vspace*{3pt}}
\mathrm{O} \bigl(\tau^2 n^{-1/8} \bigr) & \mbox{if } 0 <
\rho_- < 1 < \rho_+ < \infty. \end{array} %
\rright.
\]
Part (b) of Theorem~\ref{mprisk} follows immediately.

To prove Theorem~\ref{mprisk}(a) we show that, in fact, $\Delta
= \mathrm{O}(n^{-1/2})$ if $0 < \rho_- < \rho_+ < 1$ or $1 < \rho_- < \rho_+ <
\infty$. First, suppose that $0 < \rho_- < \rho_+ < 1$. Then, for $0
< c < (1 - \sqrt{d/n})^2$,
\[
m_{d/n} \biggl(-\frac{d}{n\tau^2} \biggr) = \frac{1}{c + d/(n\tau
^2)} - \int
_c^{\infty} \frac{1}{\{s + d/(n\tau^2)\}^2} \bigl\{1 -
F_{d/n}(s) \bigr\} \,\mathrm{d}s
\]
and
\begin{eqnarray*}
\frac{n}{d}\Delta& \leq& E \biggl\{\int_0^c
\frac{1}{s + d/(n\tau
^2)} \,\mathrm{d}\F_{n,d}(s) \biggr\} + \frac{1}{c + d/(n\tau^2)} E \bigl\{
\F_{n,d}(c) \bigr\}
\\
&& {}+ \biggl\llvert\int_c^{\infty} \frac{1}{\{s + d/(n\tau^2)\}^2}
\bigl[E \bigl\{\F_{n,d}(s) \bigr\} - F_{d/n}(s) \bigr] \,\mathrm{d}s
\biggr\rrvert
\\
& \leq& E \biggl\{\int_0^c s^{-1} \,\mathrm{d}\F_{n,d}(s) \biggr\} + \frac{1}{c
+ d/(n\tau^2)} E \bigl\{\F_{n,d}(c)
\bigr\}
\\
&&{} + \frac{1}{c + d/(n\tau^2)}\sup_{s \geq c} \bigl\llvert E \bigl\{ \F
_{n,d}(s) \bigr\} - F_{d/n}(s) \bigr\rrvert
\\
& \leq& E \bigl[s_d^{-1} \1_{\{ s_d < c\}} \bigr] +
\frac{1}{c +
d/(n\tau^2)}P(s_d < c)
\\
&&{} + \frac{1}{c + d/(n\tau^2)}\sup_{s \geq c} \bigl\llvert E \bigl\{ \F
_{n,d}(s) \bigr\} - F_{d/n}(s) \bigr\rrvert
\\
& \leq& \bigl\{E \bigl(s_d^{-2} \bigr) \bigr
\}^{1/2}P(s_d < c)^{1/2} + c^{-1}P(s_d
< c) + c^{-1}\sup_{s \geq c} \bigl\llvert E \bigl\{
\F_{n,d}(s) \bigr\} - F_{d/n}(s) \bigr\rrvert,
\end{eqnarray*}
where $ s_d \geq0$ is the smallest eigenvalue of $n^{-1}X^TX$,
$\1_{D}$ is the indicator function of the event $D$, and $P(\cdot)$
denotes the probability measure induced by the joint distribution of
$(X,\ee)$. We bound the first two terms and the last term on
right-hand side above separately. Bounding the first two terms relies
on a result of Davidson and Szarek \cite
{davidson2001local}. Their Theorem II.13, which
is a consequence of concentration of measure, implies that

%
\begin{equation}
\label{lemma2a} P(s_d \leq c) \leq\exp \biggl\{-\frac{n(1 - \sqrt {d/n})^2}{2}
\biggl(1 - \frac{c^{1/2}}{1 - \sqrt{d/n}} \biggr)^2 \biggr\},
\end{equation}
provided $c \leq1 - \sqrt{d/n}$. Additionally, Lemma~\ref{c2} in
Appendix \hyperref[appC]{C} implies that $E(s_d^{-2}) = \mathrm{O}(1)$ if $n - d > 5$. Taking $c
= (1 - \sqrt{d/n})^2/2$, it follows that
\[
\bigl\{E \bigl(s_d^{-2} \bigr) \bigr\}^{1/2}P(s_d
< c)^{1/2} + c^{-1}P(s_d < c) = \mathrm{O}
\bigl(n^{-1/2} \bigr)
\]
(in fact, we can conclude that the quantities on the left above decay
exponentially, but this is not required for the current result). It
now follows from Theorem~1.1 of Bai \textit{et al.} \cite
{bai2003convergence} that
$\Delta= \mathrm{O}(n^{-1/2})$. For the case where $1 < \rho_- < \rho_+ <
\infty$, we note that the same argument as above may be applied,
except that both $\F_{n,d}(s)$ and $F_{d/n}(s)$ have a mass of weight
$(d - n)/d$ at 0, which cancel. Theorem~\ref{mprisk}(a) follows.
\end{pf*}

\begin{pf*}{Proof of Proposition~\ref{r0inf}} Proposition~\ref{r0inf}(b)
follows directly from Proposition~\ref{rrisk}. Part (a) follows from
two applications of Jensen's inequality. If
$d +1< n $, then
\begin{eqnarray*}
R \bigl\{\hat{\bb}_r\bigl(\llVert \bb\rrVert \bigr),\bb \bigr\} & =
& E \bigl[\tr \bigl\{ \bigl(X^TX + d/\llVert \bb\rrVert
^2I_d \bigr)^{-1} \bigr\} \bigr]
\\
& \geq& d \biggl[\frac{1}{d}E \bigl\{\tr \bigl(X^TX \bigr)
\bigr\} +\frac{d}{\llVert  \bb
\rrVert  ^2} \biggr]^{-1}
\\
& = & \frac{\llVert  \bb\rrVert  ^2d/n}{\llVert  \bb\rrVert  ^2 + d/n}
\\
& = & R_r^0\bigl(\llVert \bb\rrVert ,d/n\bigr)
\end{eqnarray*}
and, since $E[\tr\{(X^TX)^{-1}\}] = d/(n-d-1)$ (Problem 3.6 of
Muirhead \cite{muirhead1982aspects}),
\begin{eqnarray*}
R \bigl\{\hat{\bb}_r\bigl(\llVert \bb\rrVert \bigr),\bb \bigr\} & =
& E \bigl[\tr \bigl\{ \bigl(X^TX + d/\llVert \bb\rrVert
^2I_d \bigr)^{-1} \bigr\} \bigr]
\\
& \leq& \frac{E[\tr\{(X^TX)^{-1}\}]}{1 + \sklfrac{1}{\llVert  \bb\rrVert  ^2}E[\tr
\{(X^TX)^{-1}\}]}
\\
& = & \frac{\llVert  \bb\rrVert  ^2d/(n-d-1)}{\llVert  \bb\rrVert  ^2 + d/(n-d-1)}
\\
& = & R_r^0 \bigl\{\llVert \bb\rrVert ,d/(n-d-1) \bigr
\}.
\end{eqnarray*}
Thus, $R_r^0\{\llVert  \bb\rrVert  ,d/(n-d-1)\} \leq R\{\hat{\bb}_r(\llVert  \bb\rrVert  ),\bb
\} \leq R_r^0(\llVert  \bb\rrVert  ,d/n)$.
It follows that if $d/n \to0$, then
\[
\sup_{\bb\in\R^d} \biggl\llvert\frac{R\{\hat{\bb}_r(\llVert  \bb
\rrVert  ),\bb\}
}{R_r^0(\llVert  \bb\rrVert  ,d/n)} - 1 \biggr\rrvert
\to0.
\]
\upqed
\end{pf*}

\begin{pf*}{Proof of Proposition~\ref{mrequiv}} Suppose that
$\hat{\bb} = \hat{\bb}(\y,X) \in\E$ and that $\bb\in
S^{d-1}(\tau)$. Let $\mathbf{e}_1 =
(1,0,\ldots,0) \in\R^d$ denote the first standard basis vector and let
$U \in\O(d)$ satisfy $\bb= \tau U\mathbf{e}_1$. Then, since
$\hat{\bb} \in\E$ and $(X,\ee)$ has the same
distribution as $(XU,\ee)$,
\begin{eqnarray*}
R(\hat{\bb},\bb) & =& E_{\bb} \bigl(\llVert \hat{\bb} - \bb\rrVert
^2 \bigr)
\\
& =& E_{\bb} \bigl(\bigl\llVert U^T\hat{\bb}(\y,X) - \tau
\mathbf {e}_1\bigr\rrVert ^2 \bigr)
\\
& =& E_{\bb} \bigl(\bigl\llVert \hat{\bb}(\y,XU) - \tau
\mathbf{e}_1\bigr\rrVert ^2 \bigr)
\\
& =& E \bigl(\bigl\llVert \hat{\bb}(XU\tau\mathbf{e}_1 + \ee,XU) -
\tau \mathbf{e}_1\bigr\rrVert ^2 \bigr)
\\
& =& E \bigl(\bigl\llVert \hat{\bb}(X\tau\mathbf{e}_1 + \ee,X) -
\tau \mathbf{e}_1\bigr\rrVert ^2 \bigr)
\\
& =& E_{\tau\mathbf{e}_1} \bigl(\llVert \hat{\bb} - \tau\mathbf{e}_1
\rrVert ^2 \bigr)
\\
& =& R(\hat{\bb},\tau\mathbf{e}_1).
\end{eqnarray*}
Part (a) of the proposition follows.

To prove part (b), we first show that\vspace*{1.5pt}
%
%
\begin{equation}
\label{prop3a} r(\tau) = \inf_{\hat{\bb} \in\E} \sup_{\bb\in
S^{d-1}(\tau)}
R(\hat{\bb},\bb).
\end{equation}
Given an estimator $\hat{\bb}$ (not necessarily orthogonally
equivariant), define\vspace*{1.5pt}
\[
\hat{\bb}_{\O}(\y,X) = \int_{\O(d)} U\hat{\bb}(
\y,XU) \,\mathrm{d}\pi_{\O(d)}(U),
\]
where $\pi_{\O(d)}$ is the uniform (Haar) measure on $\O(d)$.
Then $\hat{\bb} \in\E$ and, since $X$ and $XU$ have the same
distribution for any $U\vspace*{1pt}
\in\O(d)$,
\begin{eqnarray*}
\sup_{\bb\in S^{d-1}(\tau)} R(\hat{\bb}_{\O},\bb) & =& \sup
_{\bb
\in S^{d-1}(\tau)} E_{\bb} \biggl\{  \biggl\llVert
\int_{\O(d)} U\hat{\bb}(\y,XU) \,\mathrm{d}\pi_{\O(d)}(U) - \bb
 \biggr\rrVert^2 \biggr\}
\\
& \leq&\int_{\O(d)} \sup_{\bb\in S^{d-1}(\tau)} E \bigl\{\bigl
\llVert U \hat{\bb}(X\bb+ \ee,XU) - \bb\bigr\rrVert ^2 \bigr\} \,\mathrm{d}\pi_{\O(d)}(U)
\\
& =& \int_{\O(d)} \sup_{\bb\in S^{d-1}(\tau)} E \bigl\{\bigl
\llVert \hat{ \bb} \bigl(XU^T\bb+ \ee,X \bigr) - U^T\bb
\bigr\rrVert ^2 \bigr\} \,\mathrm{d}\pi_{\O(d)}(U)
\\
& \leq&\sup_{\bb\in S^{d-1}(\tau)} R(\hat{\bb},\bb).
\end{eqnarray*}
The identity (\ref{prop3a}) follows. Thus, by part (a) and the fact
that $\hat{\bb}_{S^{d-1}(\tau)} \in\E$,\vspace*{1pt}
\begin{eqnarray*}
r(\tau) & =& \inf_{\hat{\bb} \in\E} \sup_{\bb\in S^{d-1}(\tau)} R(\hat{
\bb},\bb)
\\
& =& \inf_{\hat{\bb} \in\E} R(\hat{\bb},\bb_1)
\\
& = &\inf_{\hat{\bb} \in\E} E_{\pi_{S^{d-1}(\tau)}} \bigl(\llVert \hat {\bb} -
\bb\rrVert ^2 \bigr)
\\
& =& E_{\pi_{S^{d-1}(\tau)}} \bigl\{\llVert \hat{\bb}_{S^{d-1}(\tau)} - \bb\rrVert
^2 \bigr\}
\\
& =& r_B(\tau),
\end{eqnarray*}
which completes the proof of the proposition.
\end{pf*}

\begin{pf*}{Proof of Theorem~\ref{adapt}} Suppose that $\bb\in
S^{d-1}(\tau)$. It is clear that (\ref{adapt2}) follows from (\ref
{adapt1}) and
Theorem~\ref{mprisk}. To prove (\ref{adapt1}), consider the risk
decomposition of the oracle and adaptive ridge
estimators
\begin{eqnarray*}
R \bigl\{\hat{\bb}_r(\tau),\bb \bigr\} & =& \biggl(\frac{d}{n}
\biggr)^2E \biggl\{ \biggl\llVert \biggl(
\frac{\tau
^2}{n}X^TX + \frac{d}{n}I_d
\biggr)^{-1}\bb  \biggr\rrVert^2 \biggr\}
\\
&&{} + \frac{1}{n^2}E \biggl\{  \biggl\llVert\tau^2
\biggl(\frac{\tau
^2}{n}X^TX + \frac{d}{n}I_d
\biggr)^{-1}X^T\ee  \biggr\rrVert^2
\biggr\},
\\
R(\check{\bb}_r,\bb) & =& \biggl(\frac{d}{n}
\biggr)^2E_{\bb} \biggl\{  \biggl\llVert \biggl(
\frac{\hat{\tau}^2}{n}X^TX + \frac{d}{n}I_d
\biggr)^{-1}\bb  \biggr\rrVert^2 \biggr\}
\\
&&{} - 2\frac{d}{n^2}E_{\bb} \biggl\{\hat{\tau}^2
\ee^TX \biggl(\frac
{\hat{\tau}^2}{n}X^TX +
\frac{d}{n}I_d \biggr)^{-2}\bb \biggr\}
\\
&&{} + \frac{1}{n^2}E_{\bb} \biggl\{  \biggl\llVert\hat{
\tau}^2 \biggl(\frac{\hat{\tau}^2}{n}X^TX +
\frac{d}{n}I_d \biggr)^{-1}X^T\ee \biggr
\rrVert ^2 \biggr\}.
\end{eqnarray*}
The triangle inequality implies
%
%
\begin{equation}
\label{adapta} \bigl\llvert R \bigl\{\hat{\bb}_r(\tau),\bb \bigr\}
- R(\check{\bb}_r,\tau) \bigr\rrvert\leq\bigl|E_{\bb}(H_1)\bigr|
+ \bigl|E_{\bb}(H_2)\bigr| + 2\bigl|E_{\bb}(H_3)\bigr|,
\end{equation}
where
\begin{eqnarray*}
H_1 & =& \biggl(\frac{d}{n} \biggr)^2 \biggl\{
 \biggl\llVert \biggl(\frac
{\tau^2}{n}X^TX +
\frac{d}{n}I_d \biggr)^{-1}\bb  \biggr
\rrVert^2 - \biggl\llVert \biggl(\frac{\hat{\tau}^2}{n}X^TX
+ \frac{d}{n}I_d \biggr)^{-1}\bb  \biggr
\rrVert^2 \biggr\},
\\
H_2 & =& \frac{1}{n^2} \biggl\{  \biggl\llVert
\tau^2 \biggl(\frac{\tau
^2}{n}X^TX +
\frac{d}{n}I_d \biggr)^{-1}X^T\ee \biggr
\rrVert ^2 -  \biggl\llVert\hat{
\tau}^2 \biggl(\frac{\hat{\tau}^2}{n}X^TX +
\frac{d}{n}I_d \biggr)^{-1}X^T\ee \biggr
\rrVert ^2 \biggr\},
\\
H_3 & =& \frac{d}{n^2}\hat{\tau}^2
\ee^TX \biggl(\frac{\hat{\tau
}^2}{n}X^TX +
\frac{d}{n}I_d \biggr)^{-2}\bb.
\end{eqnarray*}
To prove the theorem, we bound the terms $|E_{\bb}(H_1)|$, $|E_{\bb
}(H_2)|$, and
$|E_{\bb}(H_3)|$ separately.

Let $s_1 \geq\cdots\geq s_d \geq0$ denote the ordered eigenvalues
of $n^{-1}X^TX$ and let $U \in\O(d)$ be a $d \times d$ orthogonal
matrix such that $S = n^{-1}U^TX^TXU$ is diagonal. Additionally, let
$\tilde{\bb} = (\tilde{\beta}_1,\ldots,\tilde{\beta}_d)^T =
U^T\bb$
and let
$\tilde{\dd} = (\tilde{\d}_1,\ldots,\tilde{\d}_d)^T =
U^T(X^TX)^{-1/2}X^T\ee$, where $(X^TX)^{-1/2}$ denotes the
Moore--Penrose pseudoinverse of $(X^TX)^{1/2}$ if $X^TX$ is not
invertible. Then
\begin{eqnarray*}
|H_1| & =& \biggl(\frac{d}{n} \biggr)^2 \Biggl
\llvert\sum_{j = 1}^d \biggl\{
\frac{\tilde{\beta}_j^2}{(\hat{\tau}^2s_j + d/n)^2} - \frac
{\tilde{\beta}_j^2}{(\tau^2s_j + d/n)^2} \biggr\} \Biggr\rrvert
\\
& =& \biggl(\frac{d}{n} \biggr)^2 \Biggl\llvert\sum
_{j = 1}^d \frac{\tilde{\beta}_j^2s_j(\tau^2 -
\hat{\tau}^2)}{(\hat{\tau}^2s_j + d/n) (\tau^2s_j
+ d/n)} \biggl(
\frac{1}{\hat{\tau}^2 s_j + d/n} + \frac{1}{\tau^2 s_j + d/n} \biggr) \Biggr\rrvert.
\end{eqnarray*}
Since $(ax + b)^{-1} \leq(a + b)^{-1}\max\{x^{-1}, 1\}$ for
$a,b,x \geq0$,\vspace*{1pt}
\begin{eqnarray*}
|H_1| & \leq& \biggl(\frac{d}{n} \biggr)^2\sum
_{j = 1}^{d \wedge n} \biggl\{\frac{\tilde{\beta}_j^2\llvert\tau
^2 -
\hat{\tau}^2\rrvert}{(\hat{\tau}^2 + d/n) (\tau^2
+ d/n)}
\biggl(\frac{1}{\hat{\tau}^2 + d/n} + \frac{1}{\tau^2 + d/n} \biggr) \biggl(\frac{1}{s_j^2} +
s_j \biggr) \biggr\}
\\[1pt]
& \leq& \biggl(\frac{d}{n} \biggr)^2\frac{\llvert\tau^2 -
\hat{\tau}^2\rrvert}{\hat{\tau}^2 + d/n} \biggl(
\frac{1}{\hat
{\tau}^2 + d/n} + \frac{1}{\tau^2 + d/n} \biggr) \biggl(\frac
{1}{s_{d\wedge n}^2} +
s_1 \biggr).
\end{eqnarray*}
Similarly, we have\vspace*{1pt}
\begin{eqnarray*}
|H_2| & =& \frac{1}{n} \Biggl\llvert\sum
_{j = 1}^d \biggl\{\frac{\hat{\tau}^4s_j\tilde{\d}_j^2}{(\hat
{\tau}^2s_j +
d/n)^2} -
\frac{\tau^4s_j\tilde{\d}_j^2}{(\tau^2s_j +
d/n)^2} \biggr\} \Biggr\rrvert
\\[1pt]
& =& \frac{1}{n} \Biggl\llvert\sum_{j = 1}^d
\frac{(d/n)\tilde{\d}_j^2s_j(\hat{\tau}^2 - \tau^2)}{(\hat{\tau
}^2s_j +
d/n)(\tau^2s_j + d/n)} \biggl(\frac{\hat{\tau}^2}{\hat{\tau}^2s_j
+ d/n} + \frac{\tau^2}{\tau^2s_j + d/n} \biggr) \Biggr
\rrvert
\\[1pt]
& \leq&\frac{1}{n}\sum_{j = 1}^{d\wedge n}
\frac{(d/n)\tilde{\d}_j^2\llvert\hat{\tau}^2 -
\tau^2\rrvert}{(\hat{\tau}^2 + d/n)(\tau^2 + d/n)} \biggl(\frac
{1}{s_j} + s_j \biggr)
\\[1pt]
& \leq&\frac{d}{n^2}\llVert \tilde{\dd}\rrVert ^2
\frac{\llvert\hat{\tau}^2 -
\tau^2\rrvert}{(\hat{\tau}^2 + d/n)(\tau^2 +
d/n)} \biggl( \frac{1}{s_{d \wedge n}} + s_1 \biggr).
\end{eqnarray*}
Repeated application of H\"{o}lder's inequality and Lemmas \ref{c2},
\ref{c3} and \ref{c5} (found in Appendix~\hyperref[appC]{C}) imply that\vspace*{1pt}
%
%
\begin{equation}
\label{adaptb} \bigl|E_{\bb}(H_1)\bigr| + \bigl|E_{\bb}(H_2)\bigr|
= \mathrm{O} \biggl(\frac{1}{\tau^2 +
1}n^{-1/2} \biggr).
\end{equation}

To bound $|E_{\bb}(H_3)|$, we condition on $X$ and use integration by
parts (Stein's lemma,
e.g., Lemma~3.6 of Tsybakov \cite
{tsybakov2009introduction}):\vspace*{1pt}
\begin{eqnarray*}
E_{\bb}(H_3) & =& \frac{d}{n^2} E_{\bb}
\biggl\{\hat{\tau}^2 \ee^TX \biggl(\frac{\hat{\tau}^2}{n}X^TX
+ \frac{d}{n}I_d \biggr)^{-2}\bb \biggr\}
\\[1pt]
& =& \frac{2d}{n^3} E_{\bb} \biggl[\y^TX \biggl(
\frac{d}{n}I_d - \frac
{\hat{\tau}^2}{n}X^TX \biggr)
\biggl(\frac{\hat{\tau}^2}{n}X^TX + \frac{d}{n}I_d
\biggr)^{-3} \bb\1_{\{\llVert  \y\rrVert  ^2 \geq n\}} \biggr]
\\
& =& \frac{2d}{n^3} E_{\bb} \Biggl[\sum
_{j = 1}^d \frac{(ns_j\tilde{\beta}_j + n^{1/2}s_j^{1/2}\tilde{\d
}_j)(d/n -
\hat{\tau}^2s_j)\tilde{\beta}_j}{(\hat{\tau}^2s_j + d/n)^3} \1_{\{
\llVert  \y\rrVert  ^2 \geq n\}}
\Biggr].
\end{eqnarray*}
It follows that
%
%
\begin{eqnarray}\label{adaptc}
\bigl|E_{\bb}(H_3)\bigr| & \leq&\frac{2d}{n^3}E_{\bb}
\Biggl\{\sum_{j = 1}^d \biggl\llvert
\frac{(ns_j\tilde{\beta}_j +
n^{1/2}s_j^{1/2}\tilde{\d}_j)\tilde{\beta}_j}{(\hat{\tau}^2s_j
+ d/n)^2} \biggr\rrvert \Biggr\}\nonumber
\\
& \leq&\frac{2d}{n^2}E_{\bb} \Biggl\{\sum
_{j = 1}^d \frac{s_j\tilde{\beta}_j^2}{(\hat{\tau}^2s_j
+ d/n)^2} \Biggr\} +
\frac{2d}{n^{5/2}}E_{\bb} \Biggl\{\sum_{j =
1}^d
\biggl\llvert\frac{s_j^{1/2}\tilde{\d}_j\tilde{\beta}_j}{(\hat
{\tau}^2s_j
+ d/n)^2} \biggr\rrvert \Biggr\}
\\
 & =& \mathrm{O} \biggl(\frac{1}{\tau^2 + 1}n^{-1} \biggr),\nonumber
\end{eqnarray}
where we have used Lemmas \ref{c2} and \ref{c3} to obtain the last bound.
The theorem follows from (\ref{adapta}) and (\ref{adaptc}).
\end{pf*}

\section*{Appendix B}\label{appB}

This appendix is devoted to a proof of Theorem~\ref{main}, which is
fairly involved. Our first step is
to show that the minimax problem (\ref{mr}) may be reformulated as a minimax
problem for an equivalent sequence model. Ultimately, this will
substantially simplify notation and allow for a direct application of
results from Marchand \cite
{marchand1993estimation} that are important for
Theorem~\ref{main}.

\setcounter{section}{2}
\setcounter{subsection}{0}
\subsection{An equivalent sequence model}\label{sec4.2.1}

Let $\S$ be
a random orthogonally invariant $m \times m$ positive semidefinite matrix
with rank $m$, almost surely (by orthogonally invariant, we mean that
$\S$ and $U\S U^T$ have the same distribution for any $U \in\O(m)$).
Additionally, let $\dd\sim N(0,I_m)$ be an $m$-dimensional Gaussian
random vector that is independent of $\S$. Suppose that the observed
data are $(\w,\S)$, where
%
%
\begin{equation}
\label{seq0} \w= \bth+ \S^{1/2}\dd\in\R^m
\end{equation}
and $\bth\in\R^m$ is an unknown parameter.

For an estimator
$\hat{\bth} = \hat{\bth}(\w,\S)$, define the risk under squared
error loss
\[
R_{\mathrm{seq}}(\hat{\bth},\bth) = E_{\bth}\llVert \hat{\bth} - \bth
\rrVert ^2,
\]
where, abusing notation, the expectation $E_{\bth}(\cdot)$ is taken
with respect to $(\dd,\S)$ and the subscript $\bth$ indicates that
$\w= \bth+ \S^{1/2}\dd$ (we will sometimes drop the subscript
$\bth$ in $E_{\bth}(\cdot)$ if the integrand does not depend on
$\bth$). To distinguish $E_{\bth}(\cdot)$ from expectations $E_{\bb
}(\cdot)$
considered elsewhere in the paper, we emphasize that all expectations
considered in this section (Appendix \hyperref[appB]{B}) refer to the sequence
model (\ref{seq0}).

\subsection{Equivalence with the linear model}\label{sec4.2.2}

Most of the key concepts initially introduced in the context of the
linear model (\ref{lm}) have
analogues in the sequence model (\ref{seq0}).
Define
\[
\hat{\bth}_{S^{m-1}(\tau)} = E_{\pi_{S^{m-1}(\tau)}}(\bth|\w,\S),
\]
to be the posterior mean of $\bth$ under the assumption that $\bth
\sim\pi_{S^{m-1}(\tau)}$ is uniformly distributed on $S^{m-1}(\tau)$
and define
\[
\hat{\bth}_r(\tau) = E_{N(0,\tau^2/mI_m)}(\bth|\w,\S) = \tau
^2/m \bigl(\S+ \tau^2/mI_m
\bigr)^{-1}\w
\]
to be the posterior mean under the assumption that $\bth\sim
N(0,\tau^2/mI_m)$ (for both of these Bayes estimators we assume that
$\bth$ is independent of
$\dd$ and $\S$). The estimators $\hat{\bth}_{S^{m-1}(\tau)}(\tau
)$ and
$\hat{\bth}_r(\tau)$ are analogous to the
minimax estimator $\hat{\bb}_{S^{d-1}(\tau)}$ and the optimal ridge estimator
$\hat{\bb}_r(\tau)$ in the linear model, respectively. Now define the
minimax risk over $S^{m-1}(\tau)$ for the sequence model
\[
r_{\mathrm{seq}}(\tau) = \inf_{\hat{\bth}} \sup
_{\bth\in S(\tau)} R_{\mathrm{seq}}(\hat{\bth},\bth),
\]
where the infimum is over all measurable estimators for $\bth$. We
have the following analogue to Proposition~\ref{mrequiv}(b).

%
\begin{lemma}\label{b0} Suppose that $\tau\geq0$ and that $\bth
_1,\bth_2 \in
S^{m-1}(\tau)$. Then $R_{\mathrm{seq}}\{\hat{\bth}_{S^{m-1}(\tau)},\bth_1\} =
R_{\mathrm{seq}}\{\hat{\bth}_{S^{m-1}(\tau)},\bth_2\}$ and
\[
r_{\mathrm{seq}}(\tau) = \sup_{\bth\in S^{m-1}(\tau)} R_{\mathrm{seq}}\{\hat{
\bth}_{S^{m-1}(\tau)},\bth\}.
\]
\end{lemma}

The proof of Lemma~\ref{b0} is essentially the same as that of
Proposition~\ref{mrequiv} and is omitted. The next result gives an
equivalence between the linear model (\ref{lm}) and the sequence model
(\ref{seq0}) when $d \leq n$.

%
\begin{lemma} \label{lemmab1} Suppose that $m=d \leq n$ and that $\S=
(X^TX)^{-1}$.
Let $\tau\geq0$. If $\bth, \bb\in S^{m-1}(\tau)$, then
\[
R_{\mathrm{seq}} \bigl\{\hat{\bth}_r(\tau),\bth \bigr\} = R \bigl
\{ \hat{\bb}_r(\tau),\bb \bigr\}
\]
and
\[
r_{\mathrm{seq}}(\tau) = R_{\mathrm{seq}}\{\hat{\bth}_{S^{m-1}(\tau)},\bth\} =
R\{\hat{\bb}_{S^{d-1}(\tau)},\bb\} = r(\tau).
\]
\end{lemma}

Lemma~\ref{lemmab1} follows directly upon identifying $\w$ with
$\hat{\bb}_{\mathrm{ols}} = (X^TX)^{-1}X^T\y= \bb+ \linebreak[4] (X^TX)^{-1}X^T\ee$.
Lemma~\ref{lemmab1} implies that
it suffices to consider the sequence model (\ref{seq0}) (in
particular, $R_{\mathrm{seq}}\{\hat{\bth}_r(\tau),\bth\}$ and
$R_{\mathrm{seq}}\{\hat{\bth}_{S^{m-1}(\tau)},\bth\}$) in order to prove
Theorem~\ref{main}(a). Note
that Lemma~\ref{lemmab1} does not apply when $d >n$. Indeed, if $d >
n$, then the usual OLS estimator is not defined (moreover, if one uses a
pseudoinverse in place of $(X^TX)^{-1}$, then $(X^TX)^{-1}X^TX\bb$ is
not necessarily in $S^{d-1}(\tau)$). The case where $d > n$ is
considered separately below.

\subsection{Proof of Theorem \texorpdfstring{\protect\ref{main}(a)}{2(a)}}\label{sec4.2.3}

In
this section, we prove Theorem~\ref{main}(a) by bounding
%
%
\begin{equation}
\label{boundthis} R_{\mathrm{seq}} \bigl\{\hat{\bth}_r(\tau),\bth
\bigr\} - R_{\mathrm{seq}}\{\hat{\bth}_{S^{m-1}(\tau)},\bth\}.
\end{equation}
By Lemma~\ref{lemmab1}, this is equivalent to bounding
$R\{\hat{\bb}_r(\tau),\bb\} - R\{\hat{\bb}_{S^{d-1}(\tau)}(\tau
),\bb\}$.
The lower bound
%
%
\begin{equation}
\label{lboundthis} 0 \leq R_{\mathrm{seq}} \bigl\{\hat{\bth}_r(\tau),
\bth \bigr\} - R_{\mathrm{seq}}\{\hat{\bth}_{S^{m-1}(\tau)},\bth\}
\end{equation}
follows immediately from Lemma~\ref{b0}. Marchand \cite{marchand1993estimation}
obtained an upper bound on (\ref{boundthis}) in
the case where $\S= \nu^2I_m$ for fixed $\nu^2 > 0$ (i.e., in the
Gaussian sequence model with i.i.d. errors), which is one of the keys
to the proof of
Theorem~\ref{main}(a).

%
\begin{lemma}[(Theorem~3.1 from
Marchand \cite
{marchand1993estimation})]\label{lemmab2}
Suppose that $\S= \nu^2I_m$ for some fixed $\nu^2 >0$ and that $\bth
\in S^{m-1}(\tau)$. Then
\[
R_{\mathrm{seq}} \bigl\{\hat{\bth}_r(\tau),\bth \bigr
\}-R_{\mathrm{seq}}\{\hat{\bth}_{S^{m-1}(\tau)},\bth\} \leq\frac{1}{m}
\frac{\tau^2\nu^2m}{\tau^2 + \nu^2m} = \frac
{1}{m}R_{\mathrm{seq}} \bigl\{\hat{
\bth}_r(\tau),\bth \bigr\}.
\]
\end{lemma}

Thus, in the Gaussian sequence model with i.i.d. errors, the risk of
$\hat{\bth}_r(\tau)$ is
nearly as small as that of $\hat{\bth}_{S^{m-1}(\tau)}$.
Marchand's result relies on somewhat delicate calculations involving
modified Bessel functions (Robert \cite
{robert1990modified}). A direct
approach to bounding (\ref{boundthis}) for general $\S$
might involve attempting to mimic these calculations.
However, this seems daunting (Bickel \cite
{bickel1981minimax}).
Brown's identity, which relates the risk of a Bayes
estimator to the Fisher information, allows us to sidestep these
calculations and
apply Marchand's result directly.

Define the Fisher information of a random vector $\bxi\in\R^m$, with
density $f_{\bxi}$ (with respect to Lebesgue measure on $\R^m$) by
\[
I(\bxi) = \int_{\R^m} \frac{\nabla f_{\bxi}(\t)\nabla
f_{\bxi}(\t)^T}{f_{\bxi}(\t)} \,\mathrm{d}\t,
\]
where $\nabla f_{\bxi}(\t)$ is the gradient of $f_{\bxi}(\t)$.
Brown's
identity has typically been used for univariate problems or
problems in the sequence model with i.i.d. Gaussian errors
(Bickel \cite{bickel1981minimax}, Brown and Gajek \cite{brown1990information},
Brown and Low \cite{brown1991information}, DasGupta \cite{dasgupta2010false}). The
next proposition is a straightforward
generalization to the correlated multivariate Gaussian setting. Its
proof is based on Stein's lemma.

%
\begin{lemma}[(Brown's identity)] \label{lemmab3} Let $I_{\S}(\bth+ \S
^{1/2}\dd)$ denote
the Fisher information of $\bth+ \S^{1/2}\dd$, conditional on $\S$, under
the assumption that $\bth\sim\pi_{S^{m-1}(\tau)}$ is independent of
$\dd$ and
$\S$. If $\bth\in S^{m-1}(\tau)$, then
\[
R_{\mathrm{seq}}\{\hat{\bth}_{S^{m-1}(\tau)},\bth\} = E \bigl\{\tr(\S ) \bigr\}
- E \bigl[\tr \bigl\{\S^2 I_{\S} \bigl(\bth+
\S^{1/2}\dd \bigr) \bigr\} \bigr].
\]
\end{lemma}

\begin{pf} Suppose that $\bth\in S^{m-1}(\tau)$ and let
\[
f(\w) = \int_{S^{m-1}(\tau)} (2\uppi)^{-m/2}\det \bigl(
\S^{-1/2} \bigr)\mathrm{e}^{-\sfrac{1}{2}(\w-
\bth)^T\S^{-1}(\w- \bth)} \,\mathrm{d}\pi_{S^{m-1}(\tau)}(\bth)
\]
be the density of $\w= \bth+ \S^{1/2}\dd$, conditional on $\S$ and
under the assumption that
$\bth\sim\pi_{S^{m-1}(\tau)}$. Then
\[
\hat{\bth}_{S^{m-1}(\tau)} = E_{\pi_{S^{m-1}(\tau)}}(\bth|\w,\S ) = \w+
\frac{\S\nabla f(\w)}{f(\w)}.
\]
It follows that
%
%
\begin{eqnarray}\label{lemmab2a}
R_{\mathrm{seq}}\{\hat{\bth}_{S^{m-1}(\tau)},\bth\} & =&
E_{\pi_{S^{m-1}(\tau)}} \bigl\{\llVert \hat{\bth}_{S^{m-1}(\tau)} - \bth\rrVert
^2 \bigr\}\nonumber
\\
& = &E_{\pi_{S^{m-1}(\tau)}} \biggl\{  \biggl\llVert\S ^{1/2}\dd+
\frac{\S
\nabla
f(\w)}{f(\w)}  \biggr\rrVert^2 \biggr\}
\nonumber
\\
& =& E \bigl\{\tr(\S) \bigr\} + 2E_{\pi_{S^{m-1}(\tau)}} \biggl\{ \frac{\dd^T\S
^{3/2}\nabla
f(\w)}{f(\w)}
\biggr\}
\nonumber
\\[-8pt]\\[-8pt]
&&{} + E_{\pi_{S^{m-1}(\tau)}} \biggl\{\frac{\nabla f(\w)^T\S
^2\nabla
f(\w)}{f(\w)^2} \biggr\}
\nonumber
\\
& =& E \bigl\{\tr(\S) \bigr\} + 2E_{\pi_{S^{m-1}(\tau)}} \biggl\{ \frac{\dd^T\S
^{3/2}\nabla
f(\w)}{f(\w)}
\biggr\}
\nonumber
\\
&&{}  +E \bigl[\tr \bigl\{\S^2 I_{\S} \bigl(
\bth+ \S^{1/2}\dd \bigr) \bigr\} \bigr].\nonumber
\end{eqnarray}
By Stein's lemma (Lemma~3.6 of Tsybakov
\cite{tsybakov2009introduction}),
%
%
\begin{eqnarray}\label{lemmab2b}
E_{\pi_{S^{m-1}(\tau)}} \biggl\{\frac{\dd^T\S^{3/2}\nabla
f(\w)}{f(\w)} \biggr\} & =&
E_{\pi_{S^{m-1}(\tau)}} \bigl[\tr \bigl\{\S^2 \nabla^2 \log f(
\w) \bigr\} \bigr]\nonumber
\\[-8pt]\\[-8pt]
& =& -E \bigl[\tr \bigl\{\S^2 I_{\S} \bigl(\bth+
\S^{1/2}\dd \bigr) \bigr\} \bigr]. \nonumber
\end{eqnarray}
Brown's identity follows by combining (\ref{lemmab2a}) and
(\ref{lemmab2b}).
\end{pf}

Using Brown's identity, Fisher information bounds may be converted to risk
bounds, and vice-versa. Its usefulness in the present context springs
from two observations: (i) The decomposition
%
%
\begin{equation}
\label{decomp} \w= \bth+ \S^{1/2}\dd= \bigl\{\bth+ (\g\sigma
_m)^{1/2}\dd_1 \bigr\} + (\S- \gamma\sigma
_mI_m)^{1/2}\dd_2,
\end{equation}
where $\dd_1,\dd_2 \iidsim N(0,I_m)$ are independent of $\S$,
$\sigma_m$ is the smallest
eigenvalue of $\S$, and $0 < \g< 1$ is a constant; and (ii) Stam's
inequality for the Fisher information of sums of
independent random variables.

%
\begin{lemma}[(Stam's inequality; this version due to
Zamir \cite{zamir1998proof})] \label{lemmab4}
Let $\u,\vv\in\R^m$ be
independent random variables that are absolutely continuous
with respect to Lebesgue measure on $\R^m$. Then
\[
\tr \bigl[\Psi I(\u+ \vv) \bigr] \leq\tr \bigl[\Psi \bigl\{I(\u )^{-1}
+ I(\vv)^{-1} \bigr\}^{-1} \bigr]
\]
for all $m \times m$ positive definite matrices $\Psi$.
\end{lemma}

Notice in (\ref{decomp}) that $\bth+ (\g\sigma_m)^{1/2}\dd_1$ may
be viewed as
an observation from the Gaussian sequence model with i.i.d. errors,
conditional on $\S$. The necessary bound on (\ref{boundthis}) is
obtained by piecing together Brown's identity, the
decomposition (\ref{decomp}), and Stam's inequality, so that
Marchand's inequality (Lemma~\ref{lemmab2}) may be applied to $\bth+
(\g\sigma_m)^{1/2}\dd_1$.

%
\begin{lemma}\label{lemmab5} Suppose that $\S$ has rank $m$ with probability
$1$ and that $\bth\in S^{m-1}(\tau)$. Let $\sigma_1 \geq\cdots\geq
\sigma_m \geq0$ denote the eigenvalues of
$\S$. Then
\[
R_{\mathrm{seq}} \bigl\{\hat{\bth}_r(\tau),\bth \bigr\} -
R_{\mathrm{seq}}\{\hat{\bth}_{S^{m-1}(\tau)},\bth\} \leq E \biggl[ \biggl(
\frac{\sigma_1}{m\sigma_m} \wedge1 \biggr)\tr \bigl\{ \bigl(\S ^{-1} + m/
\tau^2I_m \bigr)^{-1} \bigr\} \biggr].
\]
\end{lemma}

\begin{pf} Since $\S$ is orthogonally invariant and independent of
$\dd$,
%
%
\begin{eqnarray}\label{normalRisk}
R_{\mathrm{seq}} \bigl\{\hat{\bth}_r(\tau),\bth \bigr\}
& =& E_{\bth} \bigl\{  \bigl\llVert\tau^2/m
\bigl(\S+ \tau^2/mI_m \bigr)^{-1}\w- \bth \bigr
\rrVert ^2 \bigr\}\nonumber
\\
& =& E \bigl\{  \bigl\llVert\S \bigl(\S+
\tau^2/mI_m \bigr)^{-1}\bth  \bigr
\rrVert^2 \bigr\}\nonumber
\\
&&{} + E \bigl\{ \bigl\llVert\tau^2/m \bigl(\S+
\tau^2/mI_m \bigr)^{-1}\S^{1/2}\dd
 \bigr\rrVert^2 \bigr\}\nonumber
\\
& =& E \bigl[\tr \bigl\{\tau^2/m\S^2 \bigl(\S+
\tau^2/mI_m \bigr)^{-2} \bigr\} \bigr]
\\
&&{} +E \bigl[\tr \bigl\{ \bigl(\tau^2/m \bigr)^2
\S \bigl(\S+ \tau^2/mI_m \bigr)^{-2} \bigr\}
\bigr]\nonumber
\\
& =& E \bigl[\tr \bigl\{\tau^2/m\S \bigl(\S+
\tau^2/mI_m \bigr)^{-1} \bigr\} \bigr]\nonumber
\\
 & =& E \bigl[\tr \bigl\{ \bigl(\S^{-1} + m/
\tau^2I_m \bigr)^{-1} \bigr\} \bigr].\nonumber
\end{eqnarray}
Thus, Brown's identity and (\ref{normalRisk}) imply
\begin{eqnarray*}
R_{\mathrm{seq}} \bigl\{\hat{\bth}_r(\tau),\bth \bigr\} -
R_{\mathrm{seq}}\{\hat{\bth}_{S^{m-1}(\tau)},\bth\} & =& E \bigl[\tr \bigl\{
\S^2I_{\S} \bigl(\bth+ \S^{1/2}\dd \bigr) \bigr\}
\bigr]
\\
&&{} + E \bigl[\tr \bigl\{ \bigl(\S^{-1} + m/\tau^2I_m
\bigr)^{-1} \bigr\} \bigr] - E \bigl\{\tr(\S) \bigr\}
\\
& =& E \bigl[\tr \bigl\{\S^2I_{\S} \bigl(\bth+
\S^{1/2}\dd \bigr) \bigr\} \bigr]
\\
&&{} - E \bigl[\tr \bigl\{\S^2 \bigl(\S+ \tau^2/mI_m
\bigr)^{-1} \bigr\} \bigr].
\end{eqnarray*}
Taking $\u= \bth+ (\g\sigma_m)^{1/2}\dd_1$, $\vv= (\S- \gamma
\sigma_mI_m)^{1/2}\dd_2$, and $\Psi= \S^2$ in Stam's inequality,
where $\dd_1$, $\dd_2$, and $0
< \g< 1$ are given in (\ref{decomp}), one obtains
\begin{eqnarray*}
&&R_{\mathrm{seq}} \bigl\{\hat{\bth}_r(\tau),\bth \bigr\} -
R_{\mathrm{seq}}\{\hat{\bth}_{S^{m-1}(\tau)},\bth\}
\\
&&\quad  \leq E \bigl\{\tr \bigl(\S^2 \bigl[I_{\S} \bigl\{\bth+
( \g\sigma_m)^{1/2}\dd_1 \bigr
\}^{-1}+ \S- \g\sigma_m I_m
\bigr]^{-1} \bigr) \bigr\}
\\
&&\qquad {} - E \bigl[\tr \bigl\{\S^2 \bigl(\S+ \tau^2/mI_m
\bigr)^{-1} \bigr\} \bigr].
\end{eqnarray*}
By orthogonal invariance, $ I_{\S}\{\bth+ (\g\sigma_m)^{1/2}\dd
_1\} =
\zeta I_m$ for some
$\zeta\geq0$. Thus,
%
%
\begin{eqnarray}\label{zineq}
R_{\mathrm{seq}} \bigl\{\hat{\bth}_r(\tau),\bth \bigr\}
- R_{\mathrm{seq}}\{\hat{\bth}_{S^{m-1}(\tau)},\bth\} & \leq& E \biggl(\tr
\biggl[ \S^2 \biggl\{\S+ \biggl(\frac{1}{\zeta} - \g
\sigma_m \biggr) I_m \biggr\}^{-1} \biggr]
\biggr)\nonumber
\\[-8pt]\\[-8pt]
 &&{} - E \bigl[\tr \bigl\{\S^2 \bigl(\S+
\tau^2/mI_m \bigr)^{-1} \bigr\} \bigr].\nonumber
\end{eqnarray}
Next we bound $\zeta$. Conditioning on $\S$, applying Brown's
identity with
$\g\sigma_mI_m$ in place of $\S$, and applying Marchand's inequality
(Lemma~\ref{lemmab2}) with $\nu^2 = \g\sigma_m$, we obtain
\[
m\g^2\sigma_m^2\zeta= \tr \bigl[
\g^2\sigma_m^2I_{\S} \bigl\{\bth+
(\g\sigma_m)^{1/2}\dd_1 \bigr\} \bigr] \leq m
\g\sigma_m - \biggl(1 - \frac{1}{m} \biggr)
\frac{\tau^2\g\sigma_mm}{\tau^2 + \g\sigma_mm}.
\]
Dividing by $m\g^2\sigma_m^2$ above, it follows that
\[
\zeta\leq \biggl(\frac{1}{\g\sigma_m} \biggr) \frac{\g\sigma_m
+ \tau^2/m^2}{\g\sigma_m + \tau^2/m}.
\]
Further rearranging implies that
\[
\frac{1}{\zeta} - \g\sigma_m \geq(m - 1)\frac{\g
\sigma_m\tau^2}{\g\sigma_mm^2 + \tau^2}.
\]
Hence, combining this with (\ref{zineq}),
\begin{eqnarray*}
R_{\mathrm{seq}} \bigl\{\hat{\bth}_r(\tau),\bth \bigr\} -
R_{\mathrm{seq}}\{\hat{\bth}_{S^{m-1}(\tau)},\bth\} & \leq &E \biggl(\tr \biggl[
\S^2 \biggl\{\S+ (m - 1)\frac{\g
\sigma_m\tau^2}{\g\sigma_mm^2 + \tau^2}I_m \biggr
\}^{-1} \biggr] \biggr)
\\
&&{} - E \bigl[\tr \bigl\{\S^2 \bigl(\S+ \tau^2/mI_m
\bigr)^{-1} \bigr\} \bigr].
\end{eqnarray*}
Finally, taking $\g\uparrow1$ above yields
\begin{eqnarray*}
R_{\mathrm{seq}} \bigl\{\hat{\bth}_r(\tau),\bth \bigr\} -
R_{\mathrm{seq}}\{\hat{\bth}_{S^{m-1}(\tau)},\bth\} & \leq& E \biggl(\tr \biggl[
\S^2 \biggl\{\S+ (m - 1)\frac{
\sigma_m\tau^2}{\sigma_mm^2 + \tau^2}I_m \biggr
\}^{-1} \biggr] \biggr)
\\
&&{} - E \bigl[\tr \bigl\{\S^2 \bigl(\S+ \tau^2/mI_m
\bigr)^{-1} \bigr\} \bigr]
\\
& \leq& E \biggl[ \biggl(\frac{\sigma_1}{m\sigma_m} \wedge1 \biggr)\tr \bigl\{ \bigl(
\S^{-1} + m/\tau^2I_m \bigr)^{-1}
\bigr\} \biggr],
\end{eqnarray*}
where it is elementary to verify the second inequality upon
diagonalizing $\S$. This completes the proof of the lemma.
\end{pf}

Theorem~\ref{main}(a) follows immediately from
(\ref{lboundthis}) and Lemmas
\ref{lemmab1} and \ref{lemmab5}.

\subsection{Proof of Theorem \texorpdfstring{\protect\ref{main}(b)}{2(b)}}\label{sec4.2.4}

It remains to prove Theorem~\ref{main}(b), which is achieved through a
sequence of lemmas. Similar to the proof of Theorem~\ref{main}(a),
the initial steps involve reducing the problem from the linear model to the
sequence model. In the following lemma, we derive a basic property of
orthogonally equivariant estimators for $\bb$ (in the linear model)
when $d > n$.

%
\begin{lemma}\label{obigd}
Suppose $d > n$ and that $\hat{\bb} = \hat{\bb}(\y,X) \in\E$ is
an orthogonally
equivariant estimator for $\bb$ in the linear model (\ref{lm}).
Further suppose that $X = UDV^T$,
where $U \in\O(n)$, $D$ is an $n \times n$ diagonal matrix, and $V$
is an $n \times d$ matrix with orthonormal columns. Let $V_0$ be a
$(d-n) \times d$ matrix so that $(V\enskip  V_0) \in\O(d)$. Then
$V_0^T\hat{\bb} =0$.
\end{lemma}

\begin{pf}
Let $W \in\O(d - n)$ and let $V_W = VV^T + V_0WV_0^T \in\O(d)$. Then
%
%
\begin{equation}
\label{obigd1} \hat{\bb}(\y,X) = V_W\hat{\bb}(\y,XV_W)
= V_W\hat{\bb}(\y,X).
\end{equation}
Since (\ref{obigd1}) holds for all $W \in\O(d-n)$, we must have
$V_0^T\hat{\bb} = 0$.
\end{pf}

In the next lemma, we relate the minimax risk under the linear model
$r(\tau)$ to the risk under the sequence model.

%
\begin{lemma}\label{sbigd}
Suppose that $d > n$ and let $\tau^2
> 0$. In the
sequence model (\ref{seq0}), suppose that $m=n$ and $\S=
(XX^T)^{-1}$. For $\bth= (\th_1,\ldots,\th_d)^T \in\R^d$, let
$\bth
_n = (\th_1,\ldots,\th_n)^T
\in\R^n$ be the projection onto the first $n$ coordinates. Then
\[
r(\tau) \geq\int_{S^{d-1}(\tau)} r_{\mathrm{seq}}\bigl(\llVert
\bth_n\rrVert \bigr) \,\mathrm{d}\pi_{S^{n-1}(\tau)}(\bth) + \frac{d-n}{n}
\tau^2.
\]
\end{lemma}

\begin{pf}
By Proposition~\ref{mrequiv},
%
%
\begin{equation}
\label{sbigd1} r(\tau) = \inf_{\hat{\bb} \in\E} \int_{S^{d-1}(\tau)}
R(\hat{\bb},\bb) \,\mathrm{d}\pi_{S^{d-1}(\tau)}(\bb).
\end{equation}
Assume that $\hat{\bb} = \hat{\bb}(\y,X) \in\E$ and let $X =
UDV^T$ be the
decomposition in Lemma~\ref{obigd}. Additionally, let
$\hat{\bb}_n = (\hat{\beta}_1,\ldots,\hat{\beta}_n)^T \in\R^n$ be
the first $n$
coordinates of $\hat{\bb}$. Then, under the linear model (\ref{lm}),
\begin{eqnarray*}
\llVert \hat{\bb} - \bb\rrVert ^2 & =& \bigl\llVert V^T
\hat{\bb}(\y,X) - V^T\bb\bigr\rrVert ^2 + \bigl\llVert
V_0^T \bb\bigr\rrVert ^2
\\
& =&  \bigl\llVert\hat{\bb}_n \bigl\{UDV^T
\bb+ \ee,(UD \enskip 0) \bigr\} - V^T\bb  \bigr\rrVert
^2 + \bigl\llVert V_0^T\bb\bigr\rrVert
^2.
\end{eqnarray*}
Let $\bb_n = (\beta_1,\ldots,\beta_n)^T \in\R^n$. Integrating
$\bb$
over $S^{d-1}(\tau)$ with respect to the uniform measure, making the
change of
variables
\[
\bb\mapsto(V \enskip V_0) \lleft( %
\begin{array}
{c@{\quad}c} U^T & 0
\\
0 & I_{d-n} \end{array} %
\rright) \bb,
\]
and using the fact that $\hat{\bb} \in\E$,
it follows that
\begin{eqnarray*}
&&\int_{S^{d-1}(\tau)} \llVert \hat{\bb} - \bb\rrVert ^2 \,\mathrm{d}\pi_{S^{d-1}(\tau
)}(\bb)
\\
&&\quad = \int_{S^{d-1}(\tau)}\bigl \Vert \hat{\bb}_n \bigl\{UDV^T\bb+ \ee,(UD\enskip 0) \bigr\}
- V^T\bb\bigr\Vert^2 \,\mathrm{d}\pi_{S^{d-1}(\tau)}(\bb) +
\frac{d-n}{d}\tau^2
\\
&&\quad = \int_{S^{d-1}(\tau)} \bigl\Vert\hat{\bb}_n \bigl
\{UDU^T\bb_n + \ee,(UD\enskip 0) \bigr\}
 - U^T\bb_n \bigr\Vert^2 \,\mathrm{d}\pi_{S^{d-1}(\tau)}(\bb)
+ \frac{d-n}{d}\tau^2
\\
&&\quad
 = \int_{S^{d-1}(\tau)} \bigl\Vert\hat{\bb}_n \bigl
\{UDU^T\bb_n + \ee, \bigl(UDU^T \enskip 0 \bigr) \bigr
\} - \bb_n \bigr\Vert^2 \,\mathrm{d}\pi_{S^{d-1}(\tau)}(\bb) +
\frac{d-n}{d}\tau^2.
\end{eqnarray*}
Next, for $\w\in
\R^n$ and $n \times n$ positive definite matrices $\S$, define the
estimator for the sequence model
$\hat{\bth}(\w,\S) = \hat{\bb}_n\{\S^{-1/2}\w, (\S^{-1/2}\enskip 0)\}
$. Then, with $m = n$ and $\S= (XX^T)^{-1} = UD^{-2}U^T$,
\begin{eqnarray*}
&&\int_{S^{d-1}(\tau)} R(\hat{\bb},\bb) \,\mathrm{d}\pi_{S^{d-1}(\tau)}(\bb)
\\
&&\quad= \int
_{S^{d-1}(\tau)} R_{\mathrm{seq}}(\hat{\bth},\bth_n) \,\mathrm{d}\pi_{S^{d-1}(\tau)}(\bth) + \frac{d-n}{n}\tau^2.
\end{eqnarray*}
By equivariance, $R_{\mathrm{seq}}(\hat{\bth},\boldsymbol{\vartheta})$ is
constant over spheres $\boldsymbol{\vartheta} \in
S^{d-1}(\llVert  \bth_n\rrVert  )$, which implies that $R_{\mathrm{seq}}(\hat{\bth},\bth_n)
\geq r_{\mathrm{seq}}(\llVert  \bth_n\rrVert  )$. Hence,
%
%
\begin{equation}
\label{sbigd2} \int_{S^{d-1}(\tau)} R(\hat{\bb},\bb) \,\mathrm{d}\pi_{S^{d-1}(\tau)}(\bb) \geq\int_{S^{d-1}(\tau)} r_{\mathrm{seq}}\bigl(
\llVert \bth_n\rrVert \bigr) \,\mathrm{d}\pi_{S^{d-1}(\tau)}(\bth) +
\frac{d-n}{n} \tau^2.
\end{equation}
The lemma follows from (\ref{sbigd1}) and (\ref{sbigd2}).
\end{pf}

The proof of Theorem~\ref{main}(b) will follow from a calculation
involving Lemmas \ref{lemmab5} and \ref{sbigd}. The key part of this
calculation is contained in the following lemma.

%
\begin{lemma} \label{lemmab7}
Suppose that $2 < n < d$. Let $s_1 \geq\cdots\geq s_n \geq0$ denote
the nonzero eigenvalues
of $n^{-1}X^TX$. Then
\[
r(\tau) \geq E \biggl[ \biggl(1 - \frac{s_1}{ns_n} \biggr)\tr \biggl\{
\biggl(XX^T + \frac{n(d-2)}{\tau^2(n-2)}I_n \biggr)^{-1}
\biggr\} \biggr] + \frac
{d-n}{d}\tau^2.
\]
\end{lemma}

\begin{pf}
With $\bth_n \in\R^n$, $m = n$ and $\S= (XX^T)^{-1}$, Lemma~\ref
{lemmab5} and (\ref{normalRisk}) imply that\vspace*{-2pt}
%
%
\begin{eqnarray}\label{lbbigd}
r_{\mathrm{seq}}\bigl(\llVert \bth_n\rrVert \bigr) & =&
R_{\mathrm{seq}}\{\hat{ \bth}_{S^{n-1}(\llVert  \bth
_n\rrVert  )},\bth_n\}\nonumber
\\
& \geq& R_{\mathrm{seq}} \bigl\{\hat{\bth}_r\bigl(\llVert
\bth_n\rrVert \bigr), \bth_n \bigr\} - E \biggl[ \biggl(
\frac{s_1}{ns_n}\wedge1 \biggr)\tr \biggl\{ \biggl(XX^T +
\frac{n}{\llVert  \bth_n\rrVert  ^2}I_n \biggr)^{-1} \biggr\} \biggr]
\\
& =& E \biggl[ \biggl\{ \biggl(1 - \frac{s_1}{ns_n} \biggr)\vee0
\biggr\} \tr \biggl\{ \biggl(XX^T + \frac{n}{\llVert  \bth_n\rrVert  ^2}I_n
\biggr)^{-1} \biggr\} \biggr].\nonumber
\end{eqnarray}
Additionally, if $\bth\sim\pi_{S^{d-1}(\tau)}$, then $\bth= \tau
\z/\llVert  \z\rrVert  $ in
distribution, where $\z\sim N(0,I_d)$; using basic properties of the
chi-squared distribution, it follows that\vspace*{-2pt}
%
%
\begin{equation}
\label{invid} \int_{S^{d-1}(\tau)} \frac{1}{\llVert  \bth_n\rrVert  ^2} \,\mathrm{d}\pi_{S^{d-1}(\tau
)}(\bth) = \frac{d-2}{\tau^2(n-2)},
\end{equation}
where $\bth_n = (\th_1,\ldots,\th_n)^T \in\R^n$ is the projection of
$\bth= (\th_1,\ldots,\th_d)^T \in\R^d$ onto the first $n$ coordinates.
Thus, by (\ref{lbbigd}), Jensen's inequality and (\ref{invid}),\vspace*{-2pt}
\begin{eqnarray*}
&& \int_{S^{d-1}d(\tau)} r\bigl(\llVert \bth_n\rrVert \bigr)
\,\mathrm{d}\pi_{S^{d-1}(\tau)}( \bth)
\\
&&\quad  \geq E \biggl[ \biggl(1 - \frac{s_1}{ns_n} \biggr) \tr \biggl\{
\biggl(XX^T + \int_{S^{d-1}(\tau)} \frac{n}{\llVert  \bth_n\rrVert  ^2} \,\mathrm{d}\pi_{S^{d-1}(\tau
)}(\bth) I_n \biggr)^{-1} \biggr\} \biggr]
\\
&&\quad  \geq E \biggl[ \biggl(1 - \frac{s_1}{ns_n} \biggr)\tr \biggl\{
\biggl(XX^T + \frac{n(d-2)}{\tau^2(n-2)}I_n \biggr)^{-1}
\biggr\} \biggr].
\end{eqnarray*}
The lemma follows by combining the last inequality above with
Lemma~\ref{sbigd}.
\end{pf}

We now have the tools to complete the proof of Theorem~\ref{main}(b). Suppose that $d >n$ and $\bb\in S^{d-1}(\tau)$. Then\vspace*{-2pt}
\[
R \bigl\{\hat{\bb}_r(\tau),\bb \bigr\} = E \bigl[\tr \bigl\{
\bigl(XX^T + d/\tau^2I_n \bigr)^{-1}
\bigr\} \bigr] + \frac{d-n}{d}\tau^2.
\]
Since $ R\{\hat{\bb}_{S^{d-1}(\tau)},\bb\} = r(\tau)$, Lemma~\ref
{lemmab7} implies\vspace*{-2pt}
\begin{eqnarray*}
R \bigl\{\hat{\bb}_r(\tau),\bb \bigr\} - R\{\hat{\bb
}_{S^{d-1}(\tau)}, \bb\} & =& R \bigl\{\hat{\bb}_r(\tau),\bb \bigr\} -
r(\tau)
\\
& \leq& E \bigl[\tr \bigl\{ \bigl(XX^T + d/\tau^2I_n
\bigr)^{-1} \bigr\} \bigr]
\\
&&{} - E \biggl[ \biggl(1 - \frac{s_1}{ns_n} \biggr)\tr \biggl\{
\biggl(XX^T + \frac{n(d-2)}{\tau^2(n-2)}I_n \biggr)^{-1}
\biggr\} \biggr]
\\
& \leq&\frac{1}{n}E \biggl[\frac{s_1}{s_n}\tr \biggl\{
\biggl(XX^T + \frac{d}{\tau^2}I_n \biggr)^{-1}
\biggr\} \biggr]
\\
&&{} + \frac{2(d-n)}{\tau^2(n-2)}E \biggl[\tr \biggl\{ \biggl(XX^T +
\frac{d}{\tau^2}I_n \biggr)^{-2} \biggr\} \biggr].
\end{eqnarray*}
Theorem~\ref{main}(b) follows.

\section*{Appendix C}\label{appC}
\setcounter{section}{3}
\setcounter{lemma}{0}
%
\begin{lemma} \label{c1}
Let $s_d \geq0$ denote the smallest eigenvalue of $n^{-1}X^TX$.
Suppose that $a > 0$ is a positive real number and that $n-d \geq2a + 1$.
If $d = 1$, then $E(s_d^{-a}) \leq \mathrm{e}^a$. If $d \geq2$, then
%
%
\begin{equation}
\label{lemmaC1statement} E \bigl(s_d^{-a} \bigr) \leq2 \biggl\{
\frac{\uppi}{4}\sqrt{\frac
{n^5}{(d-1)(n-d)^2}}\mathrm{e}^{n + 1/2} \biggr
\}^{2a/(n-d+1)}.
\end{equation}
\end{lemma}

\begin{pf}
Suppose first
that $d = 1$. Then $n s_d \sim
\chi^2_n$ is a chi-squared random variable on $n$ degrees of
freedom. By Theorem~1 of Ke{\v{c}}ki{\'c} and Vasi{\'c} \cite
{keckic1971some}, which gives
convenient bounds on the ratio of two gamma functions,
\[
E \bigl(s_d^{-a} \bigr) = \frac{(n/2)^a\Gamma(n/2 - a)}{\Gamma(n/2)} \leq
\frac
{(n/2)^a(n/2-a)^{n/2 - a - 1}}{(n/2)^{n/2
-1}}\mathrm{e}^a \leq \mathrm{e}^a.
\]
This proves the first part of the lemma.

Now suppose that $d \geq2$. Suppose further that (\ref
{lemmaC1statement}) is true for
$a =1$. If $0 < a_0 < 1$, then
\[
E \bigl(s_d^{-a_0} \bigr) \leq \bigl\{E
\bigl(s_d^{-1} \bigr) \bigr\}^{a_0} \leq2 \biggl\{
\frac{\uppi
}{4}\sqrt{ \frac{n^5}{(d-1)(n-d)^2}}\mathrm{e}^{n + 1/2} \biggr
\}^{2a_0/(n-d+1)}
\]
and (\ref{lemmaC1statement}) holds for $a = a_0$. Thus, we may assume
that $a \geq1$. Let $t > 0$ be a fixed positive number.
Then
%
%
\begin{equation}
\label{lemmaC1a} E \bigl(s_d^{-a} \bigr) \leq E
\bigl[s_d^{-a}\1_{\{ s_d \leq t\}} \bigr] + t^{-a}.
\end{equation}
Muirhead \cite{muirhead1982aspects}
(Corollary~3.2.19) gives the joint density
of the ordered eigenvalues, $s_1 > \cdots> s_d > 0$, of $n^{-1} X^TX$:
\[
f_{d,n}(s_1,\ldots,s_d) = c_{d,n}\exp
\Biggl(-\frac{n}{2}\sum_{j = 1}^d
s_j \Biggr)\prod_{j = 1}^d
s_j^{(n - d - 1)/2}\prod_{i < j}(s_i
- s_j),
\]
where
\[
c_{d,n} = \frac{\pi^{d^2/2}}{(2/n)^{dn/2}\Gamma_d(d/2)\Gamma_d(n/2)}
\]
and
\[
\Gamma_d(n/2) = \pi^{d(d-1)/4}\prod
_{j = 1}^d\Gamma \bigl\{(n - j + 1)/2 \bigr\}
\]
is the multivariate gamma function. Let $T_d = \{(s_1,\ldots,s_d) \in
\R^d; s_1 > \cdots> s_d > 0\}$. Then,
\begin{eqnarray*}
E \bigl[s_d^{-a}\1_{\{s_d < t\}} \bigr] & =& \int
_{T_d \cap\{s_d < t\}} s_d^{-a}f_{d,n}(s_1,
\ldots,s_d) \,\mathrm{d}s_1 \cdots \,\mathrm{d}s_d
\\
& \leq&\int_{T_{d-1}} \biggl\{\int_0^t
s_d^{-a} f_{d,n}(s_1,
\ldots,s_d) \,\mathrm{d}s_d\biggr\} \,\mathrm{d}s_1 \cdots
\,\mathrm{d}s_{d-1}
\\
& \leq&\frac{c_{d,n}}{c_{d-1,n}}\int_{T_{d-1}} \Biggl(\prod
_{j = 1}^{d-1} s_j \Biggr)^{1/2}
f_{d-1,n}(s_1,\ldots,s_{d-1}) \,\mathrm{d}s_1
\cdots \,\mathrm{d}s_{d-1}
\\
&&{} \cdot\int_0^t s^{(n-d-1)/2-a}
\mathrm{e}^{-ns/2} \,\mathrm{d}s
\\
& \le&\frac{c_{d,n}}{c_{d-1,n}} E \bigl\{\det \bigl(n^{-1}Z^TZ
\bigr)^{1/2} \bigr\}t^{(n-d+1)/2 - a},
\end{eqnarray*}
where $Z$ is an $n \times(d-1)$-dimensional matrix with i.i.d. $N(0,1)$
entries and the last inequality above follows from the fact that $n-d
\geq2a + 1$.
It is easy to check that
\[
\frac{c_{n,d}}{c_{n,d-1}} = \frac{\sqrt{\uppi}(n/2)^{n/2}}{\Gamma
(d/2)\Gamma\{(n - d + 1)/2\}}.
\]
Additionally, it is well known (Problem 3.11 in Muirhead \cite{muirhead1982aspects}, for instance) that
\[
E \bigl\{\det \bigl(n^{-1}Z^TZ \bigr)^{1/2} \bigr
\} = (2/n)^{(d-1)/2}\frac{\Gamma\{(n + 1)/2\}}{\Gamma\{(n - d + 1)/2\}}.
\]
By Corollary~1.2 of Batir \cite
{batir2008inequalities} (a variant of
Stirling's approximation),
\[
x^x\mathrm{e}^{-x}\sqrt{2x +1} \leq\Gamma(x+1)\leq
x^x\mathrm{e}^{-x}\sqrt{\uppi(2x + 1)}, \qquad \mbox{for all } x \geq0.
\]
It follows that,
\begin{eqnarray*}
\frac{c_{n,d}}{c_{n,d-1}} E \bigl\{\det \bigl(n^{-1}Z^TZ
\bigr)^{1/2} \bigr\} & =& \frac{\sqrt{\uppi}(n/2)^{(n-d+1)/2}\Gamma\{
(n +
1)/2\}}{\Gamma(d/2)\Gamma\{(n - d + 1)/2\}^2 }
\\
& \le&\frac{ \uppi n^{(n-d+2)/2} (n-1)^{(n-1)/2}\mathrm{e}^{(n-d-3)/2}
}{4(d-2)^{(d-2)/2} \sqrt{d-1} (n-d-1)^{n-d-1} (n-d) }
\\
& \le&\frac{\uppi}{4}\sqrt{\frac{n^5}{(d-1)(n-d)^2}}\mathrm{e}^{n + 1/2}
\end{eqnarray*}
and
\[
E \bigl[s_d^{-a}\1_{\{s_d < t\}} \bigr] \leq
t^{(n-d+1)/2-a}\frac{\uppi}{4}\sqrt{\frac{n^5}{(d-1)(n-d)^2}}
\mathrm{e}^{n + 1/2}.
\]
Thus, by (\ref{lemmaC1a})
\[
E \bigl(s_d^{-a} \bigr) \leq t^{(n-d+1)/2-a}
\frac{\uppi}{4}\sqrt{\frac
{n^5}{(d-1)(n-d)^2}}\mathrm{e}^{n + 1/2}+
t^{-a}.
\]
Taking $t = [(\uppi/4)\sqrt{n^{5}/\{(d-1)(n-d)^{2}\}
}\mathrm{e}^{n+1/2} ]^{-2/(n-d+1)}$ gives (\ref{lemmaC1statement}).
\end{pf}

%
\begin{lemma} \label{c2}
Let $s_1 \geq s_d \geq0$ denote the largest and smallest eigenvalues
of $n^{-1}X^TX$, respectively. Suppose that $a > 0$ is a fixed
positive real number and that
$0 < d/n \leq\rho_+ < 1$ for some fixed constant $\rho_+ \in\R$.
\begin{enumerate}[(b)]
\item[(a)] $E(s_1^a) = \mathrm{O}(1)$.

\item[(b)] If $n-d > 2a + 1$, then $E(s_d^{-a}) = \mathrm{O}(1)$.
\end{enumerate}
The constants implicit in the bounds from parts \textup{(a)} and \textup{(b)}
depend on the exponent $a$.
\end{lemma}

\begin{pf} Part (a) is well known and may be easily derived from
large deviations results for $s_1$ (see, e.g., Theorem II.13 of
Davidson and Szarek \cite{davidson2001local}). Part (b)
follows directly from Lemma~\ref{c1}.
\end{pf}

%
\begin{lemma} \label{c3}
Let $a > 0$ be a fixed positive real number. If $n > 2a$, then
\[
\sup_{\bb\in S^{d-1}(\tau)} E_{\bb} \biggl\{ \biggl(
\frac
{1}{\hat
{\tau}^2 + d/n} \biggr)^a \biggr\}= \mathrm{O} \biggl\{ \biggl(
\frac{
n/d+1}{\tau^2+1} \biggr)^a \biggr\},
\]
where the implicit constant in the big-$\mathrm{O}$ bound depends on the
exponent $a$.
\end{lemma}

\begin{pf}
Suppose that $\bb\in S^{d-1}(\tau)$. Since $\llVert  \y\rrVert  ^2/(\tau^2 + 1)
\sim\chi^2_n$ has a chi-squared
distribution with $n$ degrees of freedom,
\begin{eqnarray*}
E_{\bb} \biggl\{ \biggl(\frac{1}{\hat{\tau}^2 + d/n} \biggr)^a
\biggr\} & \leq& E_{\bb} \biggl\{ \biggl(\frac{ n/d+1}{\hat{\tau
}^2 + 1}
\biggr)^a \biggr\}
\\
& \leq&( n/d+1)^an^aE_{\bb} \bigl(\llVert \y
\rrVert ^{-2a} \bigr)
\\
& =& \mathrm{O} \biggl\{ \biggl(\frac{ n/d+1}{\tau^2+1} \biggr)^a \biggr
\}.
\end{eqnarray*}
\upqed
\end{pf}

%
%
\begin{lemma} \label{c4} Let $P_{\bb}(\cdot)$ denote the probability
measure induced by the joint distribution of $(\y,X)$, where $\y=
X\bb+ \ee$. Then
\[
\sup_{\bb\in S^{d-1}(\tau)} P_{\bb} \bigl(\hat{\tau}^2= 0
\bigr) \leq \mathrm{e}^{(-\sfrac{n}{4}) (\afrac{\tau^2}{\tau^2 + 1} )^2}.
\]
\end{lemma}

\begin{pf} Suppose that $\bb\in S^{d-1}(\tau)$. Let $t \geq0$ be
fixed. Since $V = \llVert  \y\rrVert  ^2/(\tau^2 + 1) \sim\chi^2_n$ has a chi-squared
distribution with $n$ degrees of freedom, it follows that
\[
P_{\bb} \bigl(\hat{\tau}^2 = 0 \bigr) = P_{\bb}
\biggl(V \leq\frac{n}{\tau^2 + 1} \biggr) \leq \mathrm{e}^{\afrac{nt}{\tau
^2 + 1}}E_{\bb}
\bigl(\mathrm{e}^{-tV} \bigr) = \biggl(\frac{\mathrm{e}^{\afrac{2t}{\tau^2 + 1}}}{1
+ 2t}
\biggr)^{n/2}.
\]
Taking $t = \tau^2/2$ and using the fact that $(1 - x)\mathrm{e}^{x} \leq
\mathrm{e}^{-x^2/2}$ for all $x \geq0$ yields
\[
P_{\bb} \bigl(\hat{\tau}^2 = 0 \bigr) \leq \biggl(
\frac{\mathrm{e}^{\afrac{\tau^2}{\tau^2 + 1}}}{\tau^2+1} \biggr)^{n/2} \leq \mathrm{e}^{(-\sfrac{n}{4}) (\afrac{\tau^2}{\tau^2 + 1} )^2}.
\]
\upqed
\end{pf}

%
%
\begin{lemma}\label{c5}
Suppose $a > 0$ is a fixed positive real number. Then
\[
\sup_{\bb\in S^{d-1}(\tau)} E_{\bb} \bigl( \bigl\llvert\hat{
\tau}^2 - \tau^2 \bigr\rrvert^a \bigr) =
\mathrm{O} \biggl(\frac{\tau^{2a} + 1}{n^{a/2}} \biggr),
\]
where the implicit constant in the big-$\mathrm{O}$ bound depends on the
exponent $a$.
\end{lemma}

\begin{pf} Suppose that $\bb\in S^{d-1}(\tau)$. From the
definition of $\hat{\tau}^2$,
%
%
\begin{equation}
\label{lemmaC4a} E_{\bb} \bigl( \bigl\llvert\hat{\tau}^2 -
\tau^2 \bigr\rrvert^a \bigr) \leq E_{\bb} \biggl
\{ \biggl\llvert\frac{1}{n}\llVert \y\rrVert ^{2} - \bigl(
\tau^2 + 1 \bigr) \biggr\rrvert^{a} \biggr\} +
\tau^{2a} P_{\bb} \bigl(\hat{\tau}^2 = 0 \bigr).
\end{equation}
Since $\llVert  \y\rrVert  ^2/(\tau^2 + 1) \sim\chi^2_n$,
%
%
\begin{equation}
\label{lemmaC4b} E_{\bb} \biggl\{ \biggl\llvert\frac{1}{n}\llVert
\y\rrVert ^{2} - \bigl(\tau^2 + 1 \bigr) \biggr\rrvert
^{a} \biggr\} =\mathrm{O} \biggl(\frac{\tau^{2a} +
1}{n^{a/2}} \biggr).
\end{equation}
Additionally, Lemma~\ref{c4} implies
%
%
\begin{eqnarray}\label{lemmaC4c}
\tau^{2a}P_{\bb} \bigl(\hat{\tau}^2 = 0
\bigr) & \leq&\tau^{2a} \mathrm{e}^{(-\sfrac{n}{4}) (\afrac{\tau^2}{\tau^2 +
1} )^2}\nonumber
\\
& =& \bigl(\tau^2 + 1 \bigr)^a \biggl(
\frac{\tau
^2}{\tau^2 + 1} \biggr)^a \mathrm{e}^{(-\sfrac{n}{4}) (\afrac{\tau^2}{\tau
^2 +
1} )^2}
\\
& \le& \bigl(\tau^2 + 1 \bigr)^a \biggl(
\frac{2a}{n} \biggr)^{a/2} \mathrm{e}^{-a/2}.
\nonumber
\end{eqnarray}
The lemma follows by combining (\ref{lemmaC4a}) and (\ref{lemmaC4c}).
\end{pf}
\end{appendix}

\section*{Acknowledgements}
The author would like to thank Alan Edelman, Bill Strawderman, Cun-Hui Zhang
and Sihai Zhao for their helpful comments and inspiration. The author
thanks the Associate Editor and
the referees for their careful reading of the paper and their suggestions
that helped to greatly improve its presentation.
Supported by NSF Grant DMS-1208785.


%

\printhistory

\end{document}